\documentclass[11pt,reqno]{amsart} 
\usepackage{amsmath}
\usepackage{amsthm}
\usepackage{amssymb}
\usepackage{amsfonts}
\usepackage{latexsym}
\usepackage{mathrsfs}
\usepackage{multicol}
\usepackage{fancyvrb}
\usepackage{amssymb,latexsym,amsbsy,amsthm,amsxtra,amscd,amsfonts,amsthm,amsmath,verbatim}

\usepackage{lmodern}
\usepackage[T1]{fontenc}

\usepackage[margin=1in]{geometry}
\usepackage{paralist}

\usepackage{color}

\usepackage[colorlinks]{hyperref}

\makeatletter

\@addtoreset{equation}{section}
\def\theequation{\thesection.\@arabic \c@equation}

\def\@citecolor{blue}
\def\@linkcolor{blue}
\def\@urlcolor{blue}

\def\theenumi{\@alph\c@enumi}

\makeatother

\usepackage{amssymb,latexsym,color,amsxtra,amscd,amsfonts,amsthm,amsmath,verbatim,epsfig}

\usepackage{multicol}
\usepackage{fancyvrb}
\usepackage{xypic}
\usepackage[all,cmtip]{xy}
\input xy
\xyoption{all}
\usepackage[colorlinks]{hyperref}
\def\@citecolor{blue}
\def\@linkcolor{blue}
\def\@urlcolor{blue}

\def\theequation{\arabic{equation}}
\def\theequation{\thesection.\arabic{equation}}
\numberwithin{equation}{section}
\def\im{\operatorname{Im}}

\def\ker{\operatorname{ker}}

\def\depth{\operatorname{depth}}

\def\R{\mathcal R}
\def\F{\mathcal F}
\def\ZZ{\mathbb Z}
\def\Q{\mathbb Q}
\newcommand{\m}{\mathfrak m}

\newtheorem{innercustomthm}{Theorem}
\newenvironment{customthm}[1]
  {\renewcommand\theinnercustomthm{#1}\innercustomthm}
  {\endinnercustomthm}

\theoremstyle{plain}
\newtheorem{theorem}[equation]{Theorem}
\newtheorem{corollary}[equation]{Corollary}
\newtheorem{proposition}[equation]{Proposition}
\newtheorem{lemma}[equation]{Lemma}

\theoremstyle{definition}
\newtheorem{notation}[equation]{Notation}
\newtheorem{example}[equation]{Example}
\newtheorem{definition}[equation]{Definition}

\newtheorem{remark}[equation]{Remark}

\newcommand{\ncom}{\newcommand}
\ncom{\bib}{\bibitem}

\ncom{\limns}{\underset{ s \longrightarrow \infty }{\lim}}
\ncom{\limnr}{\underset{ r \longrightarrow \infty }{\lim}}
\ncom{\beqn}{\begin{eqnarray*}}
\ncom{\eeqn}{\end{eqnarray*}}
\ncom{\beq}{\begin{eqnarray}}
\ncom{\eeq}{\end{eqnarray}}
\ncom{\f}{\frac}
\ncom{\been}{\begin{enumerate}}
\ncom{\eeen}{\end{enumerate}}
\ncom{\olin}{\overline}
\ncom{\ulin}{\underline}

\ncom{\lrar}{\longrightarrow}
\begin{document}

\title[]{An extension of Rees' theorem and Two interpretations of {a vector in the Joint reduction lattice} }
 \author{Clare D'Cruz}
 \address{Chennai Mathematical Institute, Plot H1 SIPCOT IT Park, Siruseri, 
Kelambakkam 603103, Tamil Nadu, 
India
} 
\email{clare@cmi.ac.in}
 
 \author{Shreedevi K. Masuti}
\address{Department of Mathematics, Indian Institute of Technology Dharwad, WALMI Campus, PB Road, Dharwad - 580011, Karnataka, India}
\email{shreedevi@iitdh.ac.in}
\date{\today}
 \thanks{}
 \keywords{Hilbert coefficients, bigraded Hilbert coefficients, 
joint reductions, Rees algebra, local cohomology, modified Koszul homology. 
}
\thanks{Both the authors are partially supported by a grant from Infosys Foundations.  SKM is supported by INSPIRE faculty award funded by Department of Science and Technology, Govt. of India.}
\subjclass[2010]{Primary: 13D40, 13D45, Secondary: 13A30, 14B05, 13J30}

\begin{abstract}
In \cite{rees} Rees gave a characterization for the  normal joint reduction number zero of two 
$\m$-primary ideals  in an analytically unramified Cohen-Macaulay local ring of dimension two. 
 Rees' result  is a generalization of Zariski's product theorem for complete ideals in a regular local ring of dimension two. 
 The aim of this paper is to extend Rees' theorem for the ordinary powers of $\m$-primary ideals $I$ and $J$ in a Cohen-Macaulay local ring of dimension two. 
 Following Rees' approach,
 we define the modified Koszul homology modules $M^1_{r,s}(a^k,b^k)$ for a joint reduction $(a,b)$ of $I$ and $J$. 
 Under the additional assumption that the associated graded rings of $I$ and $J$ have positive depth, we obtain  
a characterization of the joint reduction number zero of $I$ and $J$ in terms of the vanishing of the module $M^1_{0,0}(a,b)$, as well as in terms of the Hilbert coefficients  and the bigraded Hilbert coefficients. 
More generally,  we introduce the joint reduction lattice and  study the vanishing of $M^1_{r,s}(a,b)$ for any $r, s \geq 0$. This  gives   a characterization for  a vector $(r,s)$ to be in the joint reduction lattice of $I$ and $J$.  
We also  give a cohomological interpretation 
of these theorems  by investigating the  local cohomology  modules of the bigraded extended Rees algebra.  This gives another characterization for a vector $(r,s)$ to be in the joint reduction lattice and  also extends a recent result of Masuti and Verma in \cite{masuti-verma} for ordinary powers of ideals. 
 \end{abstract}

\maketitle

\section{Introduction}
Throughout this paper $(R, \m)$ is a Noetherian local ring with infinite residue field. Recall that an ideal $I$ is complete if 
$I = \olin{I}:= \{ x \in R | x^n + a_1 x^{n-1} + \cdots + a_n, a_i \in I^i \mbox{ for } i = 1, \ldots,n\}$.
In 1960, Zariski showed that if $(R,\m)$ is  a regular local ring of dimension two, then   product of complete ideals is complete (\cite{zariski-samuel}). 
In  order to generalize Zariski's result,  in 1981, Rees studied complete ideals in an analytically unramified Cohen-Macaulay local ring of dimension two.  
 To state Rees' result we need to introduce some notation. It is well-known that if $R$ is an analytically unramified local ring of dimension $d$ and $I$ an $\m$-primary ideal in $R$, 
 then   there exists a polynomial  $\olin{P}_{I}(x)  \in \Q[x]$  such that   $\olin{P}_{I}(n)  = \olin{H}_{I}(n):= \ell (R/ \olin{I^n})$ for $n \gg 0$ (\cite[Theorem~1.4]{rees0}, \cite[Theorem~1.1]{rees1}). 
 Here $\ell(M)$ denotes the length of the $R$-module $M$.  
One can write
 $P_{\olin{I}}(n):=\sum_{i=0}^{d}(-1)^i \olin{e}_i(I){n+d-i-1 \choose n+d-i}$. 
Rees proved the following interesting result: 

\begin{customthm}{1} \rm(Rees' theorem) \cite[Theorem~2.5]{rees}
\label{thm:Rees}
Let $(R,\m)$ be   an analytically unramified Cohen-Macaulay local ring of dimension two and $I,J$ be $\m$-primary ideals in $R.$ Then the following statements are equivalent:
\begin{enumerate}
 \item $\olin{e}_2 (IJ) = \olin{e}_2 (I)  + \olin{e}_2 (J) $;
 \item $\olin{I^{r+1}J^{s+1}} = a\olin{I^{r} J^{s+1}} + b \olin{I^{r+1} J^s}$ for all $r,s \geq 0.$ 
\end{enumerate}
\end{customthm}

\noindent
 In \cite{rees}, Rees proved that $\olin{e}_2(I)=0$ for every $\m$-primary ideal $I$ in a regular local ring of dimension two, thus generalizing Zariski's product theorem for complete ideals. One of the main aims of this paper is to extend Theorem~\ref{thm:Rees} for the ordinary powers of ideals $I$ and $J$. We refer to Theorem~\ref{thm:Rees} as Rees' theorem throughout this paper.

 Consider the filtrations $\F:=\{I^rJ^s\}_{r,s \in \ZZ}$ and $\olin{\F}=\{\olin{I^rJ^s}\}_{r,s \in \ZZ}$, where $I$ and $J$ are $\m$-primary ideals in $R$.  An  important tool  which was used by Rees in \cite{rees}  in generalizing Zariski's product theorem was the  bigraded normal Hilbert function
 $ H_{\olin{\F}} (r,s):= \ell ( R / \olin {I^r J^s})$, where $r,s \in \ZZ$.  In the same paper Rees showed that
in an analytically unramified local ring of dimension $d$,  there exists a polynomial ${P}_{\olin{\F}}(x,y) \in \Q[x,y]$  such that  
$P_{\olin{\F}}(r,s) = H_{\olin{\F}} (r,s) $ for all $r,s \gg 0$ and can be written as 
$$
{P}_{\olin{\F}}(r,s) = \sum_{i+j \leq d} (-1)^{i+j} {e}_{(i,j)}( \olin{\F} ) {r+i-1 \choose i} {s + j-1 \choose j}.
$$
For our purpose we consider  the bigraded Hilbert function $H_{\F} (r,s) := \ell (R/ I^rJ^s)$ for all
 $r,s \in \ZZ$. 
In \cite{bhat} P. B. Bhattacharya proved that there exists a polynomial 
 ${P}_{\F}(x,y) \in \Q[x,y]$  such that
$P_{\F}(r,s) = H_{\F} (r,s)$ for all $r,s \gg 0. $ We call this polynomial  the {\it bigraded Hilbert polynomial}. This polynomial can be written as
\beq
\label{eqn:bhat-function}
P_{\F}(r,s) = \sum_{i+j \leq d} (-1)^{i+j} {e}_{(i,j)}( \F ) {r+i-1 \choose i} {s + j-1 \choose j}.
\eeq 
We call the coefficients $e_{(i,j)}(\F)$ the {\it bigraded Hilbert coefficients}. If the ideals $I$ and $J$ are clear from the context, then for simplicity we set $e_{(i,j)}:=e_{(i,j)}(\F)$.  Let   $P_I(x) \in \mathbb Q[x]$ be   the  Hilbert-Samuel polynomial of $I$, i.e.,  $P_I(x)$ is the polynomial such that $P_{I}(n) = H_I(n) :=\ell(R/I^n)$ for $n \gg 0$. This polynomial can be written as 
$$
 P_I(n)=e_0(I)\binom{n+d-1}{d}-e_1(I)\binom{n+d-2}{d-1}+\cdots+(-1)^de_d(I) 
$$
for some integers $e_i(I)$, for $i=0,\ldots,d$, known as the {\it Hilbert coefficients of $I$}. 

In \cite{rees} Rees introduced the modified Koszul homology modules $M^1_{r,s}$ for the filtration $\olin{\F}$ and gave an explicit formula for the normal Hilbert coefficients in terms of these modules. 
In order to generalize Rees' theorem we define the modified  bigraded Koszul complex $C_{\bullet}((a^k,b^k),r,s)$ (\ref{main complex}) for $a \in I$ and $b \in J$. Using this complex,  we define the modified Koszul homology modules $M^1_{r,s}(a^k,b^k)$ for all $ 
~r,s \geq 0$ and $k \geq 1,$ and elements $a\in I $ and $b \in J$ (Definition~\ref{modified first homology}).  We set  {$M^1_{r,s}:=M^1_{r,s}(a,b)$}.
Let $(R, \m)$ be a Cohen-Macaulay local ring,  $I$ and $J$ be $\m$-primary ideals in $R$ and $(a,b)$ a joint reduction of $I$ and $J$.    If either $k \gg 0$ or $r,s \gg 0$ then the asymptotic behaviour of the  
$M^1_{r,s}(a^k,b^k)$  plays an important role in  understanding the coefficients of   $P_{\F}(r,s)$.  
Recall that $(a,b)$ is a {\it joint reduction} of $I$ and $J$ if $a \in I$, $b\in J$ and 
\beq
I^{r+1}J^{s+1}=aI^rJ^{s+1}+bI^{r+1}J^s \mbox{ for some and hence for all }r,s \gg 0
\eeq
(\cite{rees3}).   
We first prove that if $r,s \geq 0$, then  $\ell( M^1_{r,s}(a^k,b^k))$ is a polynomial of degree at most one  in $k$  for all $k \gg 0$ (Proposition~\ref{computation of length}(\ref{k large})). We explicitly describe this polynomial in terms of the Hilbert coefficients and the bigraded Hilbert coefficients.
In addition, if  we choose $a$ and $b $ to be Rees  superficial elements (see Definition~\ref{Def:Rees-superficial}),  then we study the asymptotic behaviour of $\ell(M^1_{r,s}(a^k,b^k))$ for all $k \geq 1$ and  if  either $r$ or $s$ is large 
((Proposition~\ref{computation of length}(\ref{l-r-s-a})).  
 Using this result, we show that  the differences $e_{(1,0)}-e_1(I)$ and $e_{(0,1)}-e_1(J)$ can be expressed in terms of the  length of the modules $M^1_{r,s}$  for $r, s \gg 0$,  (Proposition~\ref{computation of length}(\ref{relating e10 and e1})).  Hence,  the length of the modules $ M^1_{r,s}$ help us to  measure for the difference between $e_{(1,0)}$ and $e_1(I)$.  A similar  expression gives the difference between $e_{(0,1)}$ and $e_1(J)$.

For the filtration $\overline{\F}$,  Rees showed that ${e}_{(1,0)}( \olin{\F} ) =  \olin{e}_1(I)$
and ${e}_{(0,1)}( \olin{\F} ) =  \olin{e}_1(J)$ \cite[Theorem~1.2]{rees}. This was an important  result used in Rees' proof of Theorem~\ref{thm:Rees}. 
For the filtration $\F$, in general,  $e_{(1,0)}$ (resp. $e_{(0,1)}$) need not be equal to $e_1(I)$ (resp. $e_1(J)$). In fact, in \cite{anna-clare} the first author and A.~Guerrieri 
proved that $e_{(1,0)}=e_1(I)$ and $e_{(0,1)}=e_1(J)$ in any 
Noetherian local ring of dimension $d$. We give a counter-example to their result (Example~\ref{example-not-finite}). 
In fact, $e_{(1,0)} -e_1(I)$ (resp. $e_{(0,1)} -e_1(J)$)  can be as large as possible  (Example~\ref{example-not-finite}). This was
an obstruction in generalizing Rees' theorem for the filtration $\F.$ 
We conclude that $e_{(1,0)}\geq e_1(I)$ and $e_{(0,1)}\geq e_1(J),$ and  give a criteria for the equality in terms of the vanishing of the modules $M^1_{r,s}$ (Corollary \ref{inequality between e's}).   
As a consequence we generalize Rees' theorem in Theorem~\ref{theorem:computation of length}.

\begin{definition}
For ideals $I$ and $J$ we define the {\it joint reduction lattice} of $I$ and $J$,  denoted by $\Lambda(I|J),$ as
\[
\Lambda(I|J):=\{(r,s) \in \mathbb{N}^2: I^{r+1}J^{s+1}=aI^{r} J^{s+1}+bI^{r+1}J^{s} \mbox{ for some joint reduction } (a,b) \mbox{ of $I$ and $J$} \}.
\]
\end{definition}
We remark that if $K=\{a,c,d,b\}$ is a complete reduction of $I$ and $J$,  and  $r^{1,1}_{K}(I,J)=n$ where $r^{1,1}_{K}(I,J)$ is the joint reduction number of type $(1,1)$ with respect to $K$ 
 introduced by Hyry in \cite[Definition 3.2]{hyry},  then $(n,n) \in \Lambda(I|J).$ 
Moreover,  joint reduction vectors (with respect to a joint reduction of type $(1,1)$) introduced in \cite{SV17} are also in the joint reduction lattice.  
Recall that the ideals $I$ and $J$ are said to have {\it joint reduction number zero},  denoted by $r(I|J)=0,$ if there exists a joint reduction $(a,b)$ of $I$ and $J$ such that $IJ=aJ+bI$ (c.f. \cite{verma}).  It is clear that $r(I|J)=0$ if and only if $\Lambda(I|J)=\mathbb{N}^2.$

 In  Theorem~\ref{theorem:computation of length}, for $k \gg 0$, we characterize   joint reduction number zero for the ideals $I^k$ and $J^k$  in terms of the Hilbert and bigraded Hilbert coefficients, as well as in terms of the vanishing of the modules $M^1_{0,0}(a^k,b^k)$ for $k \gg 0.$ 
If  $G(I):=\oplus_{ n\geq 0} I^n/I^{n+1}$  and $G(J):=\oplus_{n \geq 0} J^n/J^{n+1}$ have positive depth, then Theorem~\ref{theorem:computation of length}  holds true for all $k \geq 1$ (Theorem~\ref{joint reduction number zero:KH}). 
However, if $\depth G(I)=0$ or $\depth G(J)=0$, then  the joint reduction number  of  $I$ and $J$ need not be zero even if $I$ and $J$ satisfy the equivalent conditions of Theorem~\ref{theorem:computation of length}. We demonstrate this in Example~\ref{Example:kneq1}.

More generally, we study the vanishing of the modules  $M^1_{r,s}(a^k,b^k)$ for $k \gg 0$ and $r,s \geq 0$ (Theorem~\ref{arbitrary joint reduction number:KH}). 
As a consequence, we obtain a sufficient conditions in terms of the Hilbert and bigraded Hilbert coefficients   for a vector $(r,s)$ to be in the joint reduction lattice of $I^k$ and $J^k$  for $k \gg 0$  (Corollary~\ref{Cor:AribtraryJRNklarge}). 
Under the additional assumptions (\ref{ReesSFC1}) and (\ref{ReesSFC2}), we 
obtain criteria for a  vector $(r,s)$ to be in the joint reduction lattice of $I$ and $J$ (Theorem~\ref{thm:ArbitraryJRN}).

We now describe the cohomological approach for joint reduction number zero. 
Let ${\R^\prime} (\olin{\F}):= \bigoplus_{r,s \in \mathbb Z} \olin{I^rJ^s}t_1^rt_2^s $ (resp. ${\R^\prime} ({\F}):= \bigoplus_{r,s \in \mathbb Z} {I^rJ^s}t_1^rt_2^s $) be the extended 
bigraded Rees algebra of $\olin{\F}$ (resp. $\F$) where  $t_1$ and $t_2$ are indeterminate. 
In \cite{masuti-verma} the second author and Verma  gave a new approach to Rees' theorem.  They  showed that  if $(a,b)$   is a good joint reduction of 
${\olin{\F}}$, then $\ell([H^2_{(at_1,bt_2)}(\R^{\prime}(\olin{\F}))]_{(0,0)})=-\olin{e}_2(IJ)+\olin{e}_2(I)+\olin{e}_2(J) $ (\cite[Theorem 3.7]{masuti-verma}). Moreover, in \cite[Theorem~4.1]{masuti-verma} they showed that the vanishing of $[H^2_{(at_1,bt_2)}(\R^{\prime}(\olin{\F}))]_{(0,0)}$ is equivalent to the equivalent conditions of Rees' theorem.
 This gave a cohomological interpretation of Rees' theorem. In Section~\ref{section-local cohomlogy} 
we extend these results for  the filtration $\F = \{I^rJ^s\}_{r,s \geq 0}$.

 One of the important consequences of the equalities ${e}_{(1,0)}( \olin{\F} ) =  \olin{e}_1(I)$ and ${e}_{(0,1)}( \olin{\F} ) =  \olin{e}_1(J)$ is that the local cohomology modules 
$[H^2_{(at_1,bt_2)}(\R^{\prime}(\olin{\F}))]_{(r,s)}$ have finite length for all $r,s \geq 0$, see \cite[Theorem 3.7]{masuti-verma}. 
However, for the filtration $\F$ the module $[H^2_{(at_1,bt_2)}(\R^{\prime}(\F))]_{(0,0)}$ need not have finite length (Example~\ref{example-not-finite}). In Theorem~\ref{finite length criterion}  we give necessary and sufficient conditions  for $\ell_R([H^2_{(at_1,bt_2)}(\R^{\prime} ({\F}))]_{(r,s)})$  to be finite in terms of the Hilbert and bigraded Hilbert coefficients. In particular, we show that the module $[H^2_{(at_1,bt_2)}(\R^{\prime} ({\F}))]_{(0,0)} $ has finite length if and only if ${e}_{(1,0)} =  {e}_1(I)$ and ${e}_{(0,1)} =  {e}_1(J)$ which in turn is equivalent to the vanishing of the modules $M^{1}_{i,k}$ and $M^1_{k,j}$ for all $i \geq 0$ and $j \geq 0$ and $k \gg 0$ (Corollary~\ref{corollary finite length criterion}). 

We also  we give necessary and sufficient conditions for the vanishing of  the cohomology modules
 $[H^2_{(at_1,bt_2)} (\R^{\prime}(\F))]_{(r,s)}$ (Theorem~\ref{arbitrary joint reduction number}). 
 In  fact,  we give a cohomological interpretation of Theorems \ref{arbitrary joint reduction number:KH} and \ref{thm:ArbitraryJRN} in  Theorems \ref{arbitrary joint reduction number} and \ref{thm:ArbitraryJRN-LC} respectively.
 In Theorem~\ref{thm:ArbitraryJRN-LC} we give a characterization for $(r_0,s_0) \in \Lambda(I|J)$  in terms of the vanishing of 
 $[H^2_{(at_1,bt_2)} (\R^{\prime})(\F)]_{(r_0,s_0)}$.
Putting $r_0=s_0=0$ in Theorems \ref{arbitrary joint reduction number} and \ref{thm:ArbitraryJRN-LC}  we obtain a cohomological interpretation of Theorems \ref{theorem:computation of length} and \ref{joint reduction number zero:KH} in Corollaries \ref{joint reduction number zero for large powers} and \ref{joint reduction number zero},  respectively. These results extend  Rees' theorem and \cite[Theorem 3.8]{masuti-verma} for the filtration $\F.$ We remark that Rees' theorem for joint reduction number zero has been extended for arbitrary filtration in \cite[Theorem 6.6]{MSV} under certain additional assumptions. We do  not need these additional assumptions for our results. 
{In fact, we recover  \cite[Theorem 6.6]{MSV} for the filtration $\F$ (see Corollary~\ref{joint reduction number zero})}.

In Section \ref{Section:Examples}, we give an explicit example for which $e_{(0,1)} \neq e_1(J)$ and hence $[H^2_{(at_1,bt_2)}(\R^{\prime}(\F))]_{(0,0)}$ is not finite (Example \ref{example-not-finite}).  We also give an example where both $e_{(1,0)} \neq e_1(I)$ and $e_{(0,1)} \neq e_1(J)$ (Example~\ref{Example:both}).

We conclude that the lengths of the modified Koszul  homology modules as well as the local cohomology modules are a measure  for a vector  $(r,s)$ to be in the joint reduction lattice in dimension two Cohen-Macaulay local rings. 
We hope that these approaches will be useful to find a characterization for a vector to be in the joint reduction lattice in terms of the Hilbert and bigraded Hilbert coefficients in higher dimension (see \cite{MTV15} for results in dimension 3 for the normal filtration).

We refer \cite{matsumura} for all undefined terms.

\noindent{ \bf Acknowledgement:} 
The second author  thanks the Institute of Mathematical Sciences (IMSc) and Chennai Mathematical Institute (CMI) for financial support during her post doctoral studies at IMSc and CMI, respectively, during which this work has started. 
Both the authors thank J.~K.~Verma for many fruitful discussions.

\section{Modified bigraded Koszul Complex} 
\label{preliminaries}
\noindent

In \cite[Lemma 2.2]{rees} Rees introduced the modules $M^1_{r,s}$ which played an important role in relating a  vector in the  joint reduction lattice and the  Hilbert coefficients  of  the normal filtrations $\{\overline{I^n}\}_{n \in \mathbb Z}$ and $\{\overline{J^n}\}_{n \in \mathbb Z}$ \cite[Theorem~2.4]{rees}.
 In this section,   we  construct a modified bigraded Koszul complex  $C_{\bullet}((a^k,b^k),r,s)$ for all $r,s\geq 0$ and $k \geq 1$ (c.f. \cite{marleythesis} for the modified Koszul complex of $\ZZ$-graded filtrations). 
We study  certain properties of the homology modules of this complex.  
Using this complex, for all $k \geq 1$ and $r,s \geq 0$,  we define the modified Koszul homology module $M^1_{r,s}(a^k,b^k)$ (Definition~\ref{modified first homology}).  
The  modules $M^1_{r,s}(a^k,b^k)$ are used to relate the Hilbert coefficients and the bigraded Hilbert coefficients in Section \ref{Section:JRNViaHomology}, thus extending Rees' theorem  for the filtration $\mathcal{F}$. In this section, we also  give a generalization of Huneke's fundamental Lemma (Lemma~\ref{hunekes fundamental lemma}) which relates the modules $M^1_{r,s}(a^k,b^k)$  and $H_2(C_{\bullet}((a^k,b^k),r,s))$  with the bigraded Hilbert function.

\noindent 
Let $r,s  \geq 0$, $k \geq 1$, $a \in I$ and  $b \in J$. 
We have the complexes 
$C_{\bullet}((a^k),r,s)$ and $C_{\bullet}((b^k),r,s)$  given by 
\beqn
\label{complex one}
C_{\bullet}((a^k),r,s)&:& 
0                                 \lrar \f{R}{I^rJ^s}
\buildrel .a^k \over\lrar \f{R}{I^{r+k}J^s} 
  \lrar
0\mbox{ and }\\
\label{complex two}
C_{\bullet}((b^k),r,s)&:& 
0                                 \lrar  \f{R}{I^rJ^{s}} 
\buildrel .b^k \over\lrar \f{R}{I^{r}J^{s+k}}
                         \lrar
0
\eeqn
where the maps are induced by the Koszul complex $K_{\bullet}(a^{k};R)$ and
$K_{\bullet}(b^{k};R)$,  respectively. 
We also have the chain map of complexes:
$$.b^k: C_{\bullet}((a^k),r,s) \lrar C_{\bullet}((a^k),r,s+k).$$
We call the mapping cylinder of this chain map as the modified bigraded Koszul complex 
which we denote  by  $C_{\bullet}((a^k,b^k),r,s)$ (see \cite[page 175]{rotman} for mapping cylinder).  More precisely this complex is
\beq
\label{main complex}
C_{\bullet}((a^k,b^k),r,s): 
0                                 \lrar \f{R}{I^rJ^s}
\buildrel \phi_{1} \over\lrar \f{R}{I^{r+k}J^s} \bigoplus \f{R}{I^rJ^{s+k}} 
\buildrel \phi_0 \over\lrar \f{R}{I^{r+k}J^{s+k}}
                         \lrar
0
\eeq
where the maps $\phi_0$ and $\phi_1$ are induced by the Koszul complex $K_{\bullet}(a^k,b^{k};R).$
Let $H_i(C_\bullet((a^k,b^k),r,s))$ denote the $i$-th homology of the complex 
$C_{\bullet}((a^k,b^k),r,s).$ 

\begin{theorem} 
\label{computation of homology}
Let $(R, \m)$ be a Noetherian local ring of dimension two and $I,J$ be $\m-$primary ideals 
in $R$. Let $a \in I$ and $b \in J$. Then  for all $k \geq1$ and $r,s \geq 0$,  
\been
\item
\label{computation of homology a}
${\displaystyle
H_0(C_{\bullet}((a^k,b^k),r,s))
= \f{R} {I^{r+k} J^{s+k} + (a^k, b^k)};
}$

\item
\label{computation of homology b}
${\displaystyle
H_2(C_{\bullet}((a^k,b^k),r,s))
= \f{(I^{r+k} J^{s}: (a^k))  \cap    (I^{r} J^{s+k} : (b^k))} {I^{r} J^{s} }
};$

\item
\label{computation of homology c}
If $a,b$ is a regular sequence, then 
${\displaystyle
H_1(C_{\bullet}((a^k,b^k),r,s))
= \f{  (a^k, b^k)  \cap I^{r+k} J^{s+k}  } 
{a^k I^{r} J^{s+k} +b^k I^{r+k}J^s  }
}.$
\eeen
\end{theorem}
\begin{proof} (\ref{computation of homology a}) and (\ref{computation of homology b}) are easy to verify. 

(\ref{computation of homology c}):
Consider the complex
\beq
\label{Eqn:Complex2}
C^{\prime}_{\bullet}((a^k,b^k),r,s): 
0                                 \lrar \f{R}{I^rJ^s}
\buildrel \psi_{1} \over\lrar \f{R}{I^{r+k}J^s} \bigoplus \f{R}{I^rJ^{s+k}} 
\buildrel \psi_0 \over\lrar \f{(a^k, b^k)}{ a^k I^{r}J^{s+k}  + b^k I^{r+k}J^{s} }
                         \lrar
0,
\eeq
where the maps are induced by the Koszul complex $K_{\bullet}(b^k,a^k; R)$. We claim that $\ker(\psi_0) = \im (\psi_1)$. 
We write $\olin{(.) }$ for the image of an element in respective quotients. If  $ (\olin{x}, \olin{y})   \in \ker(\psi_0)$, then  $x b^k - y a^k = a^k c + b^k d$ for some
 $c \in  I^{r}J^{s+k}$ and   $d \in  I^{r+k}
J^{s}$.  
Since $a,b$ is a regular sequence, there exists $r \in R$ such that  
$x = d + r a^k$ and $y = -c + r b^k$. Hence $(\olin{x}, \olin{y}) = \psi_1(\olin{r}) \in \im (\psi_1)$. This proves the claim. Thus we 		have the following commutative diagram with exact rows
$$
\xymatrix{
0  \ar[r]
&\ker~\psi_0
     \ar[r]^(.4){}\ar[d]
&      \f{R}{I^{r+k}J^s} \bigoplus \f{R}{I^rJ^{s+k}}     
\ar[r]^(.5){\psi_0}\ar@{=}[d]
&  \f{(a^k, b^k)}{ a^k I^{r}J^{s+k}  + b^k I^{r+k}J^{s} }
 \ar[r] \ar[d]^{\eta}
     & 0 \\
0  \ar[r]
&\ker~\phi_0\ar[r]^(.3){}
&      \f{R}{I^{r+k}J^s} \bigoplus \f{R}{I^rJ^{s+k}} 
    \ar[r]^(.5){\phi_0}
&  \f{ (a^k, b^k) + I^{r+k}J^{s+k}}   {I^{r+k}J^{s+k}}
     \ar[r] 
     & 0 .\\ }
$$
By Snake lemma   we get
\beqn
H_1(C_{\bullet}((a^k,b^k),r,s))
= \f{\ker~\phi_0}{\im~\phi_1}
= \f{\ker~\phi_0}{\im~\psi_1} = \f{\ker~\phi_0}{\ker~\psi_0}=\ker~\eta
= \f{  (a^k, b^k)  \cap I^{r+k} J^{s+k}  } 
{a^k I^{r} J^{s+k} +b^k I^{r+k}J^s  }.
\eeqn
Here the second equality is true because $\phi_1=\psi_1.$
\end{proof}

Motivated by \cite[Lemma 2.2]{rees} we introduce the modules $M^1_{r,s}(a^k,b^k)$ which 
will be useful to detect  whether a vector is in the joint reduction lattice.  

\begin{definition}
\label{modified first homology}
Let   $a \in I, ~b \in J $, $r,s \geq 0 $ and $k \geq 1$.  We define the first modified homology module to be 
\beq
\label{eqn modified first homology}
M^1_{r,s}(a^k,b^k)
:= \frac{I^{r+k} J^{s+k}}
                       {a^k I^{r}J^{s+k}+ b^k I^{r+k}J^{s}} 
                       \eeq
and we set  {$M^1_{r,s}:=M^1_{r,s}(a,b)$}.
                       \end{definition}
                  
The modules  $M^1_{r,s}(a^k,b^k)$ and  $H_1(C_{\bullet}((a^k,b^k),r,s))$ are related as follows:

\begin{lemma}
\label{lemma-fhm}
                     Let $(R,\m)$ be a Noetherian local ring of dimension two and $a,b $ a regular sequence where $a \in I$ and $b \in J$. Then for all $r,s \geq 0$ and $k \geq 1,$
 \beqn
          \ell \left(      M^1_{r,s}(a^k,b^k)\right) 
=     \ell \left(  \f{  (a^k, b^k)  +  I^{r+k} J^{s+k}  }  {(a^{k},  b^{k})  } \right) 
                     + \ell \left(  H_1(C_{\bullet}((a^k,b^k),r,s))  \right).
                     \eeqn
\end{lemma}
\begin{proof} By Theorem \ref{computation of homology} we have the short exact sequence 
 \beq
 \label{ses H and L}
 0 
 \lrar   H_1(C_{\bullet}((a^k,b^k),r,s)) 
 \lrar   M^1_{r,s}(a^k,b^k)  
 \lrar  \f{  (a^k, b^k)  +  I^{r+k} J^{s+k}  }  {(a^{k},  b^{k})  } \lrar 0.
 \eeq 
 As all the modules in  (\ref{ses H and L}) are Artinian we get the result.
\end{proof}

As a consequence of Theorem \ref{computation of homology} we obtain a bigraded version of Huneke's fundamental lemma (c.f. \cite[Lemma~2.4]{huneke}).

\begin{lemma}
\label{hunekes fundamental lemma}
\rm[Huneke's fundamental lemma for two ideals]
Let $(R, \m)$ be a Cohen-Macaulay local ring of dimension two. Let $(a,b)$  be a joint reduction of $I $ and $J$.
Then for all $r,s \geq 0$ and $k \geq 1$, 
\beq
\label{equation of hunekes fundamental lemma} \nonumber
&&\ell \left(   M^1_{r,s}(a^k,b^k) \right)  -  \ell \left(  H_2(C_{\bullet}((a^k,b^k),r,s))   \right)\\\nonumber
&=& k^2 e_{(1,1)} - H_{\F}(r+k,s+k) + H_{\F}(r,s+k)  + H_{\F}(r+k,s) -H_{\F}(r,s).
\eeq
\end{lemma}
\begin{proof}  By \cite[Theorem 2.4]{rees3}  $e_{(1,1)}=\ell (R/(a,b)).$ Hence using Theorem \ref{computation of homology}\eqref{computation of homology a} and Lemma ~\ref{lemma-fhm} we get
\beq
\label{equation of hunekes fundamental lemma-1}\nonumber
&&\ell \left(     M^1_{r,s}(a^k,b^k)\right) -  \ell \left(  H_2( C_{\bullet}((a^k,b^k),r,s))  \right)\\ \nonumber
&=& \ell \left(  \f{  (a^k, b^k)  +  I^{r+k} J^{s+k}  }  {(a^{k},  b^{k})  } \right) 
                     + \ell \left(  H_1(C_{\bullet}((a^k,b^k),r,s))  \right)-  \ell \left(  H_2( C_{\bullet}((a^k,b^k),r,s))  \right) \\ \nonumber
&=& \ell\left(\f{R}{(a^k,b^k)}\right) - \ell \left( \f{R}{  (a^k, b^k)  +  I^{r+k} J^{s+k}  }   \right) 
                     + \ell \left(  H_1(C_{\bullet}((a^k,b^k),r,s))  \right)-  \ell \left(  H_2( C_{\bullet}((a^k,b^k),r,s))  \right) \\ \nonumber
&=&\ell\left(\f{R}{(a^k,b^k)}\right) -  \ell \left(  H_0(C_{\bullet}((a^k,b^k),r,s))  \right)
                     + \ell \left(  H_1(C_{\bullet}((a^k,b^k),r,s))  \right)-  \ell \left(  H_2( C_{\bullet}((a^k,b^k),r,s))  \right) \\ \nonumber
&=&  \ell \left( \f{R} 
{(a^k,b^k )  }\right) 
- \left[\ell \left(  \f {R}{  I^{r+k} J^{s+k}  } \right) 
- \ell \left(  \f {R}{  I^{r+k} J^{s}  } \right) - \ell \left(  \f {R}{  I^{r} J^{s+k}  } \right) 
+ \ell \left(  \f {R}{  I^{r} J^{s}  } \right) 
\right] \\ \nonumber
&=&k^2 e_{(1,1)} - H_{\F}(r+k,s+k) + H_{\F}(r,s+k)  + H_{\F}(r+k,s) -H_{\F}(r,s).
\eeq
\end{proof}

We now study the properties of  $H_2(  C_{\bullet}((a^k,b^k),r,s)) $. We need a few definitions for this.

\begin{definition} 
\label{Def:Rees-superficial}
\cite[Definition 2.4]{JPV07}
(1) We say $a \in I $ (resp. $b \in J$) is a {\it Rees-superficial element} for the ideals $I$ (resp.  $J$) if 
\begin{eqnarray}
\label{superficial condition1} (a) \cap I^rJ^s&=&aI^{r-1}J^s \mbox{ for }r\gg 0 \mbox{ and all }s \geq 0 \\
 (\mbox{resp. }\label{superficial condition2}(b) \cap I^rJ^s&=&bI^rJ^{s-1} \mbox{ for }s \gg 0 \mbox{ and all }r \geq  0). 
\end{eqnarray} 
(2) We call a joint reduction $(a,b)$ of $I$ and $J$ as a  {\it Rees-joint reduction} of $I$ and $J$ if  $a$ and $b$ are Rees-superficial elements for $I$ and $J$,  respectively.   
\end{definition}

Let $(R,\m)$ be a Noetherian local ring of dimension two with infinite residue field and $I,J$ be $\m$-primary ideals in $R$. Then by \cite{rees3} there exists a Rees-joint reduction $(a,b)$ of $I$ and $J$.

The Ratliff-Rush closure of an ideal was first introduced in \cite{rat-rush} and played an important role in studying the Hilbert coefficients. The Ratliff-Rush closure for the product of ideals was considered   in \cite{jayan-verma}.
For our purpose, we introduce a refinement of the Ratliff-Rush closure for two ideals, i.e., the notion of 
Ratliff-Rush closure of $\F$ with respect to a joint reduction  (Definition~\ref{definition-rr}). This notion of the Ratliff-Rush closure is important in the study the Koszul homology $H_2( C_{\bullet}((a^k,b^k),r,s))$.

\begin{definition}
\label{definition-rr}
Let $I,J$ be $\m$-primary ideals and  $(a,b)$ a joint reduction of $I$ and $J$. 
We define the Ratliff-Rush closure of $\F$ with respect to $(a,b)$ to be  $\widetilde{\F}_{(a,b)}:=\{\widetilde{\F}_{(a,b)}(I^r,J^s)\}_{r,s \in \mathbb{Z}}$,  where
\beqn
\widetilde{\F}_{(a,b)}(I^r, J^s):=\bigcup_{k \geq 1}(I^{r+k}J^s:a^k) \cap (I^rJ^{s+k}:b^k).
\eeqn
\end{definition}

\begin{remark} One can verify that 
$I^r J^s \subseteq \widetilde{\F}_{(a,b)}(I^r, J^s) \subseteq \widetilde{I^r J^s}$, where $\widetilde{I}$ denotes the Ratliff-Rush closure of $I$  introduced in \cite{rat-rush} (see \cite{jayan-verma} for the Ratliff-Rush closure of product of ideals).
\end{remark}

\begin{remark}
\label{RR}
 For all  $k \geq 1$, 
 $                   (I^{r+k}J^s:a^k) \cap (I^rJ^{s+k}:b^k) 
  \subseteq    (I^{r+k+1}J^s:a^{k+1}) \cap (I^rJ^{s+k+1}:b^{k+1})$. As $R$ is Noetherian, 
  $
  \widetilde{\F}_{(a,b)}(I^r,J^s)=(I^{r+k}J^s:a^k) \cap (I^rJ^{s+k}:b^k)$ for some and hence  for all $k \gg 0.
  $ 
\end{remark}

 \begin{lemma} 
 \label{second homology} 
 Let $(R,\m)$ be a Cohen-Macaulay  local ring of dimension two and $I$, $J$ be $\m$-primary ideals in $R$. 
 Let $(a,b)$ be a joint reduction of $I$ and $J$.  
 \begin{enumerate}
 
\item 
 \label{second homology-a} 
 {Fix   $r,s \geq 0$.  Then there exists $k \gg 0$ (which depends on $r,s$)},  such that  $H_2(  C_{\bullet}((a^k,b^k),r,s))   =\frac{\widetilde{\F}_{(a,b)}(I^r, J^s)}{I^rJ^s}
$ and hence is independent of $k$ for $k \gg 0.$ 

\item 
 \label{second homology-b} 
Suppose  $a$ and $b$ are Rees-superficial elements for $I$ and $J$. 
  If either $r \gg 0$ or $s \gg 0$, then for all $k \geq 1$,
$H_2( C_{\bullet}((a^k,b^k),r,s))  =0$. 
 \end{enumerate}
 \end{lemma}
\begin{proof}
(\ref{second homology-a}):
 By Remark \ref{RR} 
  $
  \widetilde{\F}_{(a,b)}(I^r,J^s)=(I^{r+k}J^s:a^k) \cap (I^rJ^{s+k}:b^k)$ for all $k \gg 0.
  $ 
  Hence, by Theorem~\ref{computation of homology}, for $k \gg 0,$
  $$
              H_2  (C_{\bullet}((a^k,b^k),r,s)) 
  = \frac{\widetilde{\F}_{(a,b)}(I^r, J^s)}{I^rJ^s}. 
  $$ 
  
  (\ref{second homology-b}): 
  Since $R$ is Cohen-Macaulay and  $a,b$ is a system of parameters, it is  a regular sequence. 
  Hence for all $k \geq 1$, either $I^{r+k} J^s:a^k = I^{r}J^s$ or $I^r J^{s+k}:b^k = I^r J^s$ by  (\ref{superficial condition1}) and (\ref{superficial condition2}).
This implies that $H_2(  C_{\bullet}((a^k,b^k),r,s))  =0$.  
\end{proof}

\section{Characterization of a vector in the joint reduction lattice in terms of modified  Koszul homology}
\label{Section:JRNViaHomology}

Let $I$ and $J$ be $\m$-primary ideals in a Cohen-Macaulay local ring  $(R,\m)$ of dimension two.  
In \cite{rees} Rees used the modules  $M^1_{r,s}$ for the filtration $\overline{\F}$ to characterize the  joint reduction number zero of $\olin{\F}.$ The aim of this section is to extend Rees' theorem for the filtration $\F = \{ I^r J^s\}_{r,s \geq 0}$ and to characterize the joint reduction number zero in terms of the vanishing of the modules $M^1_{0,0}(a^k,b^k)$ for $k \gg 0$.  
We also investigate the relationship between  $e_{(1,0)} - e_1(I)$ and $e_{(0,1)} - e_1(J)$ using the modules  $M^1_{r,s}$.  
More generally,  we characterize the vector $(r,s) \in \Lambda(I|J)$ in terms of the vanishing of the modules $M^1_{r,s}(a,b)$,  and also in terms of the bigraded and the Hilbert coefficients under certain additional  assumptions.

Fix $r \geq 0$. Then  $\ell( J^s/ J^s I^r )$ is a polynomial of degree one for all large $s$, i.e., there exist integers $f_{0}(I^r)$ and  $f_{1}(I^r)$ such that for all large $s$
\beqn
 \ell( J^s/ J^s I^r ) &= &f_{0}(I^r) s - f_{1}(I^r).
\eeqn
Therefore, for all $s \gg 0$ and $r \geq 0$
  \beq
\label{length-one} \nonumber
            \ell\left( \f{R}{I^rJ^s} \right) =  \ell\left( \f{R}{J^s} \right) +  \ell\left( \f{J^s}{I^rJ^s} \right)
            &=& e_0(J) \binom{s+1}{2} - e_1(J) s + e_2(J) + f_{0}(I^r) s - f_{1}(I^r)\\
 &=&    e_0(J)\binom{s+1}{2}-g_{1}(r)s+g_{2}(r), 
 \eeq
 where   
 \beq
 \label{length-one g1 and g2}
 \begin{array}{lcllcl}
 g_{1}(r) &=&  e_1(J) -  f_{0}(I^r), &\hspace{.2in} g_{2}(r) &=&  e_2(J)  - f_{1}(I^r)   \hspace{.2in} \mbox{for all $r \geq 0$}.
 \end{array}
 \eeq 
In particular,
\beq
\label{g1(0)andg2(0)}
g_1(0)=e_1(J) \mbox{ and } g_2(0)=e_2(J).
\eeq
Comparing the equations \eqref{eqn:bhat-function} and \eqref{length-one},  for $r \gg 0$ we get 
 \beq
 \label{g for large r}
 \begin{array}{lcllcl}
 g_{1}(r) &=&  e_{(0,1)} -  e_{(1,1)}r, &\hspace{.2in} g_{2}(r) &=&  e(I){r+1 \choose 2} - e_{(1,0)}r  +e_2(IJ)    \hspace{.2in} \mbox{for all } r \gg 0.
 \end{array}
\eeq
\noindent Similarly,   fix $s\geq 0$ and for all  $r \gg 0$, there exist integers  $f_{0}'(J^s)$ and  $f_{1}'(J^s)$ such that  
\[
 \ell( I^r/ I^r J^s) = f_{0}'(J^s) r - f_{1}'(J^s)\]
 and hence for all $r \gg 0$ and $s \geq 0$
   \beq
  \label{length-two} \nonumber
            \ell\left( \f{R}{I^rJ^s} \right) 
= \ell\left( \f{R}{I^r} \right) +  \ell\left( \f{I^r}{ I^r J^s} \right)
&=& e_0(I) \binom{r+1}{2} - e_1(I) r + e_2(I) + f_{0}'(J^s) r - f_{1}'(J^s)\\
 &=&    e(I)\binom{r+1}{2}-h_1(s)r+h_2(s) 
 \eeq
where   
 \beq
 \label{length-one h1 and h2}
 \begin{array}{lcllcl}
 h_{1}(s) &=&  e_1(I) -  f_{0}'(J^s), &\hspace{.2in} h_{2}(s) &=&  e_2(I)  - f_{1}'(J^s)   \hspace{.2in} \mbox{for all } s \geq 0.
  \end{array}
   \eeq
  In particular,
  \beq
  \label{h1(0)andh2(0)}
  h_1(0) = e_1(I) \mbox{ and } h_2(0) = e_2(I).
  \eeq
 Comparing the equations \eqref{eqn:bhat-function} and \eqref{length-two},  for $s \gg 0$ we get 
 \beq 
 \begin{array}{lcllcl}
 h_{1}(s) &=&  e_{(1,0)} -  e_{(1,1)}s, &\hspace{.2in} h_{2}(s) &=&  e(J){s+1 \choose 2} - e_{(0,1)}s  +e_2(IJ)    \hspace{.2in} \mbox{for all } s \gg 0.
 \end{array}
\eeq

\noindent In the next proposition we study asymptotic behaviour of the module $M^1_{r,s}(a^k,b^k)$.

\begin{proposition}
\label{computation of length}
 Let $(R,\m)$ be a Cohen-Macaulay local ring of dimension two  and $I$, $J$ be $\m$-primary ideals in $R$. Let $(a,b)$ be a joint reduction of $I$ and $J$. In addition, we assume that  $a$ and $b$ are Rees-superficial elements  for part (\ref{l-r-s-a})
 and (\ref{relating e10 and e1}). Then the following hold true.
\begin{enumerate}
 \item 
 \label{k large}
 Let $r,s \geq 0.$ Then for $k \gg 0,$  
\begin{eqnarray*}
      \ell( M^1_{r,s}(a^k,b^k))
      &=& [e_{(0,1)}-g_1(r)-e_{(1,1)}r +e_{(1,0)}-h_1(s)-e_{(1,1)}s]k\\
  && -  e_{(1,1)}rs 
  +     [e_{(0,1)}-g_1(r)]s
  +[e_{(1,0)}-h_1(s)]r\\
 &&   -e_2(IJ)+g_2(r)+h_2(s)-\ell(R/I^rJ^s)
 + 
 \ell\left(\frac{\widetilde{\F}_{(a,b)}(I^r, J^s)}{I^rJ^s}\right).
 \end{eqnarray*}
In particular,  $\ell( M^1_{r,s}(a^k,b^k))$ is a polynomial of degree at most one  in $k$. 

\item
\label{l-r-s-a} 
Let $k \geq 1.$ Then
\beq
\label{l-r-s-a-1} 
\ell(  M^1_{r,s}(a^k,b^k)) &= &k^2e_{(1,1)}  + k \left[ g_1(r+k)- g_1(r) \right] \hspace{.2in} \mbox{for $r \geq 0$ and $s  \gg 0$},\\
\label{l-r-s-a-2} 
 \ell(M^1_{r,s}(a^k,b^k))&=&  k^2 e_{(1,1)} + k \left[ h_1(s+k)- h_1(s) \right] \hspace{.2in} \mbox{for $s \geq 0$ and $r  \gg 0$}.
 \eeq
 In particular, if $r \geq 0$ and $s \gg 0$  (resp.   $s \geq 0$ and $r \gg 0$),  then $\ell( M^1_{r,s}(a^k,b^k)  )$ is independent of $s$ (resp. $r$) and the joint reduction $(a,b)$.

 \item 
 \label{relating e10 and e1} Let $k \geq 1.$ Then for all $r,s \gg 0$
 \beq
\label{relating e10 and e1 a} 
 \sum_{i=0}^{r-1}    \ell( M^1_{ik,s}(a^k,b^k)) &=&  k  \left[e_{(0,1)}- e_1(J) \right] \mbox{ and }\\
\label{relating e10 and e1 b} 
\sum_{i=0}^{s-1}    \ell(M^1_{r,ik}(a^k,b^k)) &=& k   \left[e_{(1,0)}- e_1(I) \right].
\eeq 
 \end{enumerate}
 \end{proposition}
\begin{proof} (\ref{k large}): From Lemma~\ref{hunekes fundamental lemma}  and Lemma~\ref{second homology}(\ref{second homology-a}),  for $k \gg 0,$ we have 
\begin{eqnarray*}
&& \ell ( M^1_{r,s}(a^k,b^k))\\
&=& k^2e_{(1,1)} 
-       \left[ \ell \left( \f{R}{I^{r+k}J^{s+k}} \right)
-       \ell \left( \f{R}{I^{r+k}J^s} \right)
-       \ell \left( \f{R}{I^rJ^{s+k}} \right)
 +    \ell\left( \f{R}{I^rJ^s} \right) \right]
 + \ell\left(\frac{\widetilde{\F}_{(a,b)}(I^r, J^s)}{I^rJ^s}\right)\\
 &=& k^2e_{(1,1)}  - \left[  e_{(1,1)} (r+k) (s+k) 
 -       e_{(1,0)} (r+k) - e_{(0,1)} (s+k) + e_2(IJ)  \right]\\
 &&-   \left[
    h_1(s) (r+k) - h_2(s)  
+ g_1(r) (s+k)-g_2(r) 
 +       \ell\left( \f{R}{I^rJ^s} \right)   \right]   + \ell\left(\frac{\widetilde{\F}_{(a,b)}(I^r, J^s)}{I^rJ^s}\right)\\
 &&\hspace{3.5in} \mbox{from (\ref{eqn:bhat-function}),(\ref{length-one}), (\ref{length-two}) }\\
 &=&[e_{(0,1)}-g_1(r)-e_{(1,1)}r +e_{(1,0)}-h_1(s)-e_{(1,1)}s]k
  -  e_{(1,1)}rs 
  +     [e_{(0,1)}-g_1(r)]s
  +[e_{(1,0)}-h_1(s)]r\\
 &&   -e_2(IJ)+g_2(r)+h_2(s)-\ell(R/I^rJ^s)
 + 
 \ell\left(\frac{\widetilde{\F}_{(a,b)}(I^r, J^s)}{I^rJ^s}\right).
 \end{eqnarray*}
 
\noindent (\ref{l-r-s-a}): It is  enough to prove  (\ref{l-r-s-a-1}) as the proof of  (\ref{l-r-s-a-2}) is similar.
Let $r \geq 0$ and $s \gg 0.$ As $a$ and $b$ are Rees-superficial elements,  $H_2(  C_{\bullet}((a^k,b^k),r,s))  =0$ by Lemma \ref{second homology}(\ref{second homology-b}). Hence from Lemma~\ref{hunekes fundamental lemma} and (\ref{length-one}) we have 
\begin{eqnarray*}
 \ell( M^1_{r,s}(a^k,b^k) ) 
 &=& k^2 e_{(1,1)} + k \left[g_1(r+k) -  g_1(r)\right].
 \end{eqnarray*}

\noindent (\ref{relating e10 and e1}): 
Applying (\ref{l-r-s-a-1}), for all $s \gg 0,$  we get
\beq
\label{relating g1r and e1J} \nonumber
 \sum_{i=0}^{r-1}    \ell ( M^1_{ik,s}(a^k,b^k))
&  =& rk^2 e_{(1,1)} +  k \sum_{i=0}^{r-1}[ g_1(ik+k)-g_1(ik)]\\ \nonumber
&=&  r k^2e_{(1,1)} +k  \left[  g_1(rk)  - g_1(0) \right]\\
&=& r k^2e_{(1,1)} +k  \left[  g_1(rk)   - e_1(J) \right] \hspace{.8in} \mbox{[from (\ref{g1(0)andg2(0)})]}
\eeq
Since $g_1(r)=e_{(0,1)}-e_{(1,1)}r$ for $r \gg 0$ by (\ref{g for large r}), we get
\begin{eqnarray*}
  \sum_{i=0}^{r-1}   \ell( M^1_{ik,s}(a^k,b^k)) =  k  \left[e_{(0,1)}- e_1(J) \right].
\end{eqnarray*}

\noindent Replacing $g_1(r)$ by $h_1(s)$, $e_{(0,1)}$ by $e_{(1,0)}$ in the proof of $(\ref{relating e10 and e1 a})$ we get $(\ref{relating e10 and e1 b})$.
\end{proof}

As a corollary we express   $\ell(M^1_{0,0}(a^k,b^k))$ in terms of the  bigraded Hilbert coefficients and the Hilbert coefficients. 

\begin{corollary}
\label{cor:computation of length}
 Let $I$, $J$ be $\m$-primary ideals in a Cohen-Macaulay local ring $(R,\m)$  of dimension two and $(a,b)$ a joint reduction of $I$ and $J.$ Then for all $k \gg 0,$  $\ell(M^1_{0,0}(a^k,b^k))$ is a polynomial of degree at most one in $k$ and  this polynomial can be written as
\begin{eqnarray*}
            \ell( M^1_{0,0}(a^k,b^k))
 &=& [e_{(0,1)} - e_1(J) + e_{(1,0)}-e_1(I)]k
 -         e_2(IJ) + e_2(J)+e_2(I).
 \end{eqnarray*}
 In particular,   $\ell( M^1_{0,0}(a^k,b^k))$ is independent of the joint reduction chosen. 
 \end{corollary}
\begin{proof}  For $r=s=0$,  $ \ell\left(\frac{\widetilde{\F}_{(a,b)}(I^r, J^s)}{I^rJ^s}\right)=0$. 
We have $g_1(0)=e_1(J), ~g_2(0)=e_2(J)$ by (\ref{g1(0)andg2(0)}),  and $h_1(0)=e_1(I)$ and $h_2(0)=e_2(I)$ by (\ref{h1(0)andh2(0)}).  Hence substituting $r=s=0$ in 
Proposition~\ref{computation of length}(\ref{k large}) we get the result.
\end{proof}

In Corollary~\ref{inequality between e's}, we give a formula for the difference  $e_{(1,0)} - e_1(I)$ and $e_{(0,1)}- e_1(J)$ and a criterion for the equality to hold. This gives a generalization to \cite[Theorem~1.2]{rees}.

\begin{corollary}
\label{inequality between e's}
Let $(R,\m)$ be a Cohen-Macaulay local ring of dimension two, $I,J$ be 
 $\m$-primary ideals in $R$ and $(a,b)$ a Rees-joint reduction of $I$ and $J$. Then 
 \begin{enumerate}
 \item 
 \label{FormulaForTheDifference}
 \beq
 \label{FormulaForTheDifference1}
 e_{(0,1)}- e_1(J)  & =&\sum_{i=0}^{r-1}    \ell( M^1_{i,s})  \mbox{ for } r,s \gg 0 \mbox{ and }\\
 \label{FormulaForTheDifference2}e_{(1,0)}- e_1(I) &=& \sum_{i=0}^{s-1}    \ell(M^1_{r,i}) \mbox{ for } r,s \gg 0   .
 \eeq
\item 
\label{inequality between e's-2}
  $e_{(0,1)} \geq e_1(J)$  (resp.  $e_{(1,0)} \geq e_1(I)$) and the equality holds if and only if  for all $i \geq 0$ and $s \gg 0$, $M^1_{i,s}=0$ (resp. for all $i \geq 0$ and $r \gg 0$, $M^1_{r,i}=0$).
  \item 
  \label{inequality between e's-1}
  Let $i \geq 0$. Then for all $s \gg 0$ (resp. $r \gg 0$),  $\ell(M^1_{i,s})$ (resp. $\ell(M^1_{r,i}$))  is independent of $s$ (resp.  $r$).
\end{enumerate}
\end{corollary}
\begin{proof}
 Put  $k=1$ in Proposition~\ref{computation of length}(\ref{relating e10 and e1}) to get (\ref{FormulaForTheDifference1}) and (\ref{FormulaForTheDifference2}). (\ref{inequality between e's-2}) is immediate from (\ref{FormulaForTheDifference1}) and (\ref{FormulaForTheDifference2}). By putting $k=1$ in (\ref{l-r-s-a-1}) and (\ref{l-r-s-a-2}) we get \ref{inequality between e's}(\ref{inequality between e's-1}). 
\end{proof}

The inequalities in Corollary~\ref{inequality between e's}(\ref{inequality between e's-2}) can be strict (see Examples~\ref{example-not-finite} and \ref{Example:both}). In fact, Example~\ref{example-not-finite} shows that the difference $e_{(0,1)} - e_1(J)$ can be as large as possible. 

\begin{notation}
Let $i,j\geq 0$ and $(a,b)$ a Rees-joint reduction of $I$ and $J$. Since $\ell(M^1_{i,s})$ (resp. $\ell(M^1_{r,i}$))  is independent of $s$ (resp.  $r$) for all $s \gg 0$ (resp $r \gg 0$) by Corollary~\ref{inequality between e's}(\ref{inequality between e's-1}), we  set
\beqn
M^1_{i,*} &:=& \ell(M^1_{i,s}) \hspace*{1in} \mbox{for $s \gg 0$}\\
\label{alpha-2}M^1_{*,j} &:=& \ell(M^1_{r,j}) \hspace*{1in} \mbox{for $r \gg 0$}
\eeqn 
\end{notation}

In Proposition~\ref{nonnegativity of coefficients} we generalize Corollary~\ref{inequality between e's}(\ref{inequality between e's-2}) and prove that $ e_{(0,1)} - g_1(r) - r e_{(1,1)} \geq 0$ and $e_{(1,0)} - h_1(s) - s e_{(1,1)}  \geq 0$. We also give a criteria for the equality to hold. We need the following lemma for this purpose.

\begin{lemma}
\label{Lemma:Mi*}
 For $i,j\geq 0$ and $(a,b)$ a Rees-joint reduction of $I$ and $J$, 
 \beq
\label{alpha}
M^1_{i,*} &=&    g_1(i+1)- g_1(i) + e_{(1,1)}. \\ \label{beta}
M^1_{*,j} &=&    h_1(j+1)- h_1(j) + e_{(1,1)}.
\eeq
In particular, $g_1(i+1)- g_1(i) + e_{(1,1)} \geq 0$ and $h_1(j+1)- h_1(j) + e_{(1,1)} \geq 0.$
\end{lemma}
\begin{proof}
Put $k=1$ in (\ref{l-r-s-a-1}) and (\ref{l-r-s-a-2}) to get the result.
\end{proof}

\begin{remark}
\label{remark-alpha-i}
By Lemma~\ref{Lemma:Mi*} and (\ref{g for large r}) it follows that for $i \gg 0$, $M^1_{i,*}=0$. Similarly, for $j \gg 0$, $M^1_{*,j}=0$. 
\end{remark}

\begin{proposition}
   \label{nonnegativity of coefficients}
 Let $(R,\m)$ be a Cohen-Macaulay local ring of dimension two, $I,J$ be $\m$-primary ideals in $R$ and $(a,b)$ a Rees-joint reduction of $I$ and $J$. 
 Then 
 \begin{enumerate}
   \item 
   \label{nonnegativity of coefficients a}
   For all  $r \geq 0$, $ e_{(0,1)} - g_1(r) - r e_{(1,1)}  \geq 0$ and equality holds if and only if $M^1_{i,*}=0$ for all $i \geq r$.
      \item 
      \label{nonnegativity of coefficients b}
   For all  $s \geq 0$, $e_{(1,0)} - h_1(s) - s e_{(1,1)}  \geq 0$ and equality holds if and only if $M^1_{*,j}=0$ for all $j \geq s$.
  \end{enumerate}
 \end{proposition}
\begin{proof} (\ref{nonnegativity of coefficients a})
Since $M^1_{i,*} =   g_1(i+1)- g_1(i) + e_{(1,1)}$ by (\ref{alpha}) and $g_1(0)=e_1(J)$ by (\ref{g1(0)andg2(0)}),   
$$\displaystyle{ \sum_{i=0}^{r-1} M^1_{i,*} } = g_1(r)-g_1(0)+re_{(1,1)}=g_1(r)-e_1(J)+re_{(1,1)}.$$
Hence
\beqn
        e_{(0,1)} - g_1(r) -r e_{(1,1)} 
=e_{(0,1)}- e_1(J) - \sum_{i=0}^{r-1} M^1_{i,*} 
                      \geq  e_{(0,1)}- e_1(J) - \sum_{i\geq 0} M^1_{i,*} 
                     = 0 \hspace{.3in} \mbox{(by (\ref{FormulaForTheDifference1}))}.
\eeqn
Hence it follows that  $e_{(0,1)} - g_1(r) -r e_{(1,1)}=0$ if and only if $M^1_{i,*}=0$ for all $i \geq r$. 
The proof of (\ref{nonnegativity of coefficients b}) is similar. 
\end{proof}

We observe that if we put  $r=0$ and $s=0$  in Proposition~\ref{nonnegativity of coefficients}  we obtain Corollary~\ref{inequality between e's}(\ref{inequality between e's-2}).

We are now ready to generalize Rees' theorem for the filtration $\F.$ 
This result characterizes joint reduction number zero of $I^k$ and $J^k$ for $ k\gg 0$ in terms of the Hilbert coefficients and the bigraded Hilbert coefficients. 
We need to  make an additional assumption  that $e_{(1,0)} = e_1(I)$ and $e_{(0,1)} = e_1(J).$ These assumptions hold true   for the filtration $\overline{\F}$ (c.f. \cite[Theorem 1.2]{rees}).

\begin{theorem}
\label{theorem:computation of length}
Let $(R,\m)$ be a Cohen-Macaulay local ring of dimension two  and $I$, $J$ be $\m$-primary ideals in $R$. Then the following statements are equivalent:
\been
\item
\label{theorem:computation of length a}
$r(I^k | J^k) = 0$ for all $k \gg 0$;

\item
\label{theorem:computation of length b}
 $r(I^k | J^k) = 0$ for some $k \gg 0$;
\item
\label{theorem:computation of length c}
$e_{(1,0)} = e_1(I)$, $e_{(0,1)} = e_1(J)$, and $e_2(IJ) = e_2(I) + e_2(J);$
\item \label{theorem:computation of length d}
there exists a joint reduction $(a,b)$ of $I$ and $J$ such  that 
$M^1_{0,0}(a^k,b^k)=0$ for all $k \gg 0$.
\eeen
\end{theorem}
\begin{proof}
(\ref{theorem:computation of length a}) $\implies$ (\ref{theorem:computation of length b}) is clear.

(\ref{theorem:computation of length b}) $\implies$ (\ref{theorem:computation of length c}): 
Fix $k \geq 1$ such that  $r(I^k | J^k) = 0$.  Then $r(I^{nk} | J^{nk}) = 0$ for all $n \geq 1$.  
Hence there exists a joint reduction $(a,b)$ of $I$ and $J$ such that for all $n \gg 0$
 \beq
  \label{hfl and jr0 1} \nonumber
  0 &=& M^1_{0,0}(a^{nk},b^{nk}) \\
  &=& \nonumber [e_{(0,1)} - e_1(J)  + e_{(1,0)}-e_1(I)]nk
 -         e_2(IJ) + e_2(J)+e_2(I) \hspace{.1in} \mbox{(by Corollary~\ref{cor:computation of length})}.
 \eeq
This implies that $e_2(IJ) = e_2(J)+e_2(I)$ and $e_{(0,1)} - e_1(J)  + e_{(1,0)}-e_1(I)=0$. Since $e_{(0,1)} - e_1(J) \geq 0 $ and  $e_{(1,0)}-e_1(I) \geq 0$ by Corollary~\ref{inequality between e's}, $e_{(0,1)} = e_1(J)$   and  $e_{(1,0)}=e_1(I)$. This proves (\ref{theorem:computation of length c}). 
 
  (\ref{theorem:computation of length c}) $\implies$ (\ref{theorem:computation of length d}): 
  Let $(a,b)$ be a joint reduction of $I$ and $J$.  Then  the result follows from Corollary \ref{cor:computation of length}.

(\ref{theorem:computation of length d}) $\implies$ (\ref{theorem:computation of length a}): This is clear.
\end{proof}

We remark that
the equivalent conditions of Theorem~\ref{theorem:computation of length} need not imply that the joint reduction number of $I$ and $J$ is zero. Example~\ref{Example:kneq1} illustrates this.  Before that we make the following remark.

\begin{remark}  
\label{independent of joint reduction}
From \cite[Theorem~3.2]{verma-joint} it follows that if $r(I|J)=0$ for $\m$-primary ideals $I$ and $J$ 
in a Cohen-Macaulay local ring of dimension two with infinite residue field,  
then the condition $IJ=aJ+bI$ holds for every joint reduction of $I$ and $J$. 
\end{remark}

\begin{example}
\label{Example:kneq1}
 Let $R=k[|x,y|], I=(x^4,x^3y,xy^3,y^4)$ and $\m=(x,y)$. 
 Then  $(x^4,y)$ is a 
 joint reduction of $I$ and $\m$, and  $I^2\m^2=x^8\m^2+y^2I^2$. Hence $r(I^2| \m^2) = 0$. Therefore  the equivalent conditions of Theorem~\ref{theorem:computation of length} hold true. 
 But $r(I|\m) \neq 0,$ since $x^2y^3 \in \m I \setminus x^4 \m+y I$ (c.f. Remark~\ref{independent of joint reduction}).
\end{example}

In the next theorem we show that if we make an additional assumption that 
$\depth G(I) \geq 1$ and $\depth G(J) \geq 1,$ then Theorem~\ref{theorem:computation of length} holds  true for $k=1.$

\begin{theorem}
\label{joint reduction number zero:KH}
  Let $(R,\m)$ be a Cohen-Macaulay local ring of dimension two and $I,J$ be 
$\m$-primary ideals in $R$. Assume that $\depth G(I),\depth G(J)\geq 1$. Then following statements are equivalent:	
\begin{enumerate}
\item \label{joint reduction number zero:KH1}
$r(I|J) = 0$; 

 \item \label{joint reduction number zero:KH2}
 $e_{(1,0)}=e_1(I)$, $e_{(0,1)}=e_1(J)$ and $e_2(IJ)=e_2(I)+e_2(J)$;

  \item \label{joint reduction number zero:KH3}
there exists a  joint reduction $(a,b)$ of $I$ and $J$ such that $M^1_{0,0}(a^k,b^k)=0$ for all $k\geq 1.$
\end{enumerate}
\end{theorem}
\begin{proof}
(\ref{joint reduction number zero:KH1}) $\implies$ (\ref{joint reduction number zero:KH2}): Follows from Theorem~\ref{theorem:computation of length}.

(\ref{joint reduction number zero:KH2}) $\implies$ (\ref{joint reduction number zero:KH3}):  By \cite[Lemma~1.2]{rees3} there exists a Rees-joint reduction 
 $(a,b)$ of $I$ and $J$ such that $a \in I\setminus \m I$ and $b \in J \setminus \m J$.  
 Let $a^*$ (resp. $b^*$) denotes the image of $a$ (resp. $b$) in $[G(I)]_1$ (resp. $[G(J)]_1$). 
 Since $\depth G(I)$ (resp. $\depth G(J)$) $\geq 1$, $a^*$ (resp. $b^*$) is a 
 nonzero-divisor in $G(I)$ (resp. $G(J)$). Hence  by \cite[Lemma~2.1]{huckaba-marley},
 \begin{eqnarray*}
  (a) \cap I^n=aI^{n-1} \mbox{ and }(b) \cap J^n=bJ^{n-1} \mbox{ for all }n >0.
 \end{eqnarray*}
{To complete the proof we need to show that 
$I^kJ^k=a^kJ^k+b^kI^k$ for  all 
 $k \geq 1$.  By Theorem~\ref{theorem:computation of length} there exists $N$ so that 
   $I^kJ^k=a^kJ^k+b^kI^k$ for all $k \geq N.$ We use 
  decreasing induction on $k$ to show that $I^kJ^k=a^kJ^k+b^kI^k$ for all $k \geq 1.$ 
 
 We claim that $I^{N-1}J^{N-1}=a^{N-1}J^{N-1}+b^{N-1}I^{N-1}$. 
 We first show that 
  $I^{N-1}J^N=a^{N-1}J^N+b^NI^{N-1}$.
  Let $x \in I^{N-1}J^N$. Then $ax \in I^NJ^N.$ 
 Hence there exist   $p \in J^N$ and $q \in I^{N}$ such that 
 \beq
 \label{equation for ax}
 ax=a^Np+b^N q.
 \eeq
  This implies that $q \in (a) \cap I^{N}=aI^{N-1}$.   
 Hence there exists $q^\prime \in I^{N-1}$ such that $q=aq^\prime$. Plugging in (\ref{equation for ax}) we get 
 $ax=a^{N}p+ab^N q^\prime. $ As $a$ is a regular element 
 \beqn
  x=a^{N-1}p+b^N q^\prime \in a^{N-1}J^N+b^N I^{N-1}.
 \eeqn
 Hence  $I^{N-1}J^N=a^{N-1}J^N+b^kI^{N-1}$. Repeating the above argument we get that 
  $I^{N-1}J^{N-1}=a^{N-1}J^{N-1}+b^{N-1}J^{N-1}$.
Thus by decreasing induction on $k$ we get that  $I^kJ^k=a^kJ^k+b^kI^k$ for all $k \geq 1$.  
Hence $M^1_{0,0}(a^k,b^k)=0$ for all $k\geq 1.$
} 
 
(\ref{joint reduction number zero:KH3}) $\implies$ (\ref{joint reduction number zero:KH1}):  This follows from the definition of $M^1_{0,0}(a^k,b^k)$.
\end{proof}

In the following theorem we generalize Theorem~\ref{theorem:computation of length}. As a consequence we give a sufficient condition  for the vector $(r_0,s_0)\in \Lambda(I^k|J^k)$ for $k \gg 0$ in Corollary~\ref{Cor:AribtraryJRNklarge}.

\begin{theorem}
\label{arbitrary joint reduction number:KH}
Let $(R,\m)$ be a Cohen-Macaulay local ring of dimension two and $I,J$ be 
$\m$-primary ideals in $R$. Let $(a,b)$ be a joint reduction of $I$  and $J$. Let  $r_0,s_0 \geq 0$. Then the following statements are equivalent:
\begin{enumerate}
 \item 
\label{arbitrary joint reduction number:KH1}
for $ k\gg 0$
$$I^{r_0+k}J^{s_0+k}=a^kI^{r_0}J^{s_0+k}+b^kI^{r_0+k}J^{s_0};$$
 
 \item \label{arbitrary joint reduction number:KH2}
 $e_{(0,1)}=g_1(r_0) + r_0e_{(1,1)}$, $e_{(1,0)} = h_1(s_0) + s_0 e_{(1,1)}$ and 
 $e_2(IJ)=g_2(r_0)+h_2(s_0)-\ell(R/I^{r_0}J^{s_0})+
r_0s_0e_{(1,1)}+\ell\left(\frac{\widetilde{\F}_{(a,b)}(I^{r_0}, J^{s_0})}{I^{r_0}J^{s_0}}\right)$;

\item  \label{arbitrary joint reduction number:KH3} $M^1_{r_0,s_0}(a^k,b^k)=0$ for $k \gg 0.$
\end{enumerate}
\end{theorem}
\begin{proof}
 (\ref{arbitrary joint reduction number:KH1}) $\implies$ (\ref{arbitrary joint reduction number:KH2}):  Since $I^{r_0+k}J^{s_0+k}=a^kI^{r_0}J^{s_0+k}+b^kI^{r_0+k}J^{s_0}$ for $k \gg 0,$ by definition $M^1_{r_0,s_0}(a^k,b^k)=0$ for $k \gg 0.$ Hence by Proposition~\ref{computation of length}(\ref{k large})
$$[e_{(0,1)}-g_1(r_0) - r_0e_{(1,1)}]+[e_{(1,0)} - h_1(s_0) - s_0 e_{(1,1)}]=0.$$
Since $e_{(0,1)}-g_1(r_0) - r_0e_{(1,1)} \geq 0$ and $e_{(1,0)} - h_1(s_0) - s_0 e_{(1,1)} \geq 0$ by Proposition~\ref{nonnegativity of coefficients},  
we get that $e_{(0,1)}=g_1(r_0) + r_0e_{(1,1)}$, $e_{(1,0)} = h_1(s_0) + s_0 e_{(1,1)}$.  
Now using Proposition~\ref{computation of length}(\ref{k large}) we get
 $$e_2(IJ)=g_2(r_0)+h_2(s_0)-\ell(R/I^{r_0}J^{s_0})+
r_0s_0e_{(1,1)}+\ell\left(\frac{\widetilde{\F}_{(a,b)}(I^{r_0}, J^{s_0})}{I^{r_0}J^{s_0}}\right). $$ 
 
(\ref{arbitrary joint reduction number:KH2}) $\implies$  (\ref{arbitrary joint reduction number:KH3}): By Proposition~\ref{computation of length}(\ref{k large})   
for every joint reduction $(a,b)$ of $I$ and $J$ $M^1_{r_0,s_0}(a^k,b^k) =0$ 
 for $k \gg 0$. 
 
 (\ref{arbitrary joint reduction number:KH3}) $\implies$  (\ref{arbitrary joint reduction number:KH1}): This is clear by definition of $M^1_{r_0,s_0}(a^k,b^k) $. 
\end{proof}

\begin{corollary}
\label{Cor:AribtraryJRNklarge}
Let the assumptions be as in Theorem~\ref{arbitrary joint reduction number:KH}. 
If any of the equivalent conditions of Theorem~\ref{arbitrary joint reduction number:KH} are satisfied,  then  $(r_0,s_0)\in \Lambda(I^k|J^k)$ for all $k \gg 0$.
\end{corollary}
\begin{proof}
By assumption there exists a joint reduction $(a,b)$ of $I$ and $J$ such that 
 for $ k\gg 0$
\beq 
\label{Eqn:JRNArbitray}
I^{r_0+k}J^{s_0+k}=a^kI^{r_0}J^{s_0+k}+b^kI^{r_0+k}J^{s_0}. \eeq
Then multiplying Equation~\ref{Eqn:JRNArbitray} by $I^{(k-1)r_0}J^{(k-1)s_0}$ we get
\[
I^{kr_0 + k} J^{ks_0 + k} = a^k I^{kr_0} J^{ks_0+k} + b^k I^{kr_0 + k} J^{ks_0}.
\]
Hence $(r_0,s_0) \in \Lambda(I^k|J^k)$ for $k \gg 0.$
\end{proof}

We  now give a criterion for the  a vector $(r_0,s_0)$ to be in $\Lambda(I|J)$ under the additional assumptions (\ref{ReesSFC1}) and (\ref{ReesSFC2}).  

\begin{theorem}
\label{thm:ArbitraryJRN}
 Let $(R,\m)$ be a Cohen-Macaulay local ring of dimension two and $I,J$ be 
$\m$-primary ideals in $R$.  Let  $r_0,s_0 \geq 0$. Assume that there exists a joint reduction  $(a,b)$ of $I$ and $J$ 
such that
\begin{eqnarray}
 \label{ReesSFC1}(a) \cap I^{r_0+k}J^{s_0}&=&aI^{r_0+k-1}J^{s_0} \mbox{ for }k \geq 1  \\
 \label{ReesSFC2}(b) \cap I^{r_0}J^{s_0+k}&=&bI^{r_0}J^{s_0+k-1} \mbox{ for }k \geq 1. 
\end{eqnarray}
Then the following statements are equivalent:
\begin{enumerate}
 \item \label{thm:ArbitraryJRN1} $(r_0,s_0) \in \Lambda(I|J);$

 \item \label{thm:ArbitraryJRN2} $e_{(0,1)}=g_1(r_0) + r_0e_{(1,1)}$, $e_{(1,0)} = h_1(s_0) + s_0 e_{(1,1)}$ and 
 $e_2(IJ)=g_2(r_0)+h_2(s_0)-\ell(R/I^{r_0}J^{s_0})+
r_0s_0e_{(1,1)}+\ell\left(\frac{\widetilde{\F}_{(a,b)}(I^{r_0}, J^{s_0})}{I^{r_0}J^{s_0}}\right)$;

\item \label{thm:ArbitraryJRN3} $M^1_{r_0,s_0}(a^k,b^k)=0$ for all $k \geq 1.$
\end{enumerate}
\end{theorem}
\begin{proof}
(\ref{thm:ArbitraryJRN1}) $\implies $ (\ref{thm:ArbitraryJRN2}): Let $(a_1,b_1)$ be a joint reduction of $I$ and $J$ such that for all $r \geq r_0$ and $s \geq s_0,$
 \[
 I^{r+1}J^{s+1}=a_1 I^r J^{s+1} + b_1 I^{r+1} J^s. 
\]
Then by induction on $k$  it follows that for all $k \geq 1,$ and all $r \geq r_0$ and $s \geq s_0$
\[
 I^{r+k}J^{s+k}=a_1^k I^r J^{s+k} + b_1^k I^{r+k} J^{s}.
\] 
Hence the result follows from Theorem~\ref{arbitrary joint reduction number:KH}. 
 
(\ref{thm:ArbitraryJRN2}) $\implies $ (\ref{thm:ArbitraryJRN3}): By Proposition~\ref{computation of length}(\ref{k large}) $M^1_{r_0,s_0}(a^k,b^k)=0$ for $k \gg 0.$ Therefore 
  \[
   I^{r_0+k}J^{s_0+k}=a^k I^{r_0} J^{s_0+k} + b^k I^{r_0+k} J^{s_0}, \mbox{ say for } k \geq N.
  \]
To complete the proof we need to show  $I^{r_0+k}J^{s_0+k}=a^k I^{r_0} J^{s_0+k} + b^k I^{r_0+k} J^{s_0}$ for $ k \geq 1.$ 
We  first show that 
 \[
   I^{r_0+N-1}J^{s_0+N}=a^{N-1} I^{r_0} J^{s_0+N} + b^N I^{r_0+N-1} J^{s_0}.
  \]
Let $x \in   I^{r_0+N-1}J^{s_0+N}.$ Then $a x \in I^{r_0+N}J^{s_0+N}.$ Let $ax = a^N p + b^N q$ for some $p \in I^{r_0} J^{s_0+N} $ and $q \in I^{r_0+N} J^{s_0}.$ Then $q \in (a) \cap I^{r_0+N} J^{s_0}=a I^{r_0+N-1} J^{s_0}$ by (\ref{ReesSFC1}).  Let $q = aq' $ for some $q' \in I^{r_0+N-1} J^{s_0}.$ Then 
$x=a^{N-1} p + b^{N} q' \in a^{N-1} I^{r_0} J^{s_0+N} + b^N I^{r_0+N-1} J^{s_0}. $ 

Similar argument shows that 
 \[
   I^{r_0+N}J^{s_0+N-1}=a^{N} I^{r_0} J^{s_0+N-1} + b^{N-1} I^{r_0+N} J^{s_0}.
  \]
Continuing as above we get that for all $k \geq 1$
\[
I^{r_0+k}J^{s_0+k}=a^k I^{r_0} J^{s_0+k} + b^k I^{r_0+k} J^{s_0}.
\]  
Hence $M^1_{r_0,s_0}(a^k,b^k)=0$ for all $ k\geq 1.$ 
  
(\ref{thm:ArbitraryJRN3}) $\implies $ (\ref{thm:ArbitraryJRN1}):  This follows from the definition of $M^1_{r_0,s_0}(a^k,b^k).$ 
\end{proof}

\section{Characterization of a vector in the joint reduction lattice in terms of Local Cohomology}
\label{section-local cohomlogy}
The aim of this section is to characterize a vector in the joint reduction lattice in terms of the vanishing of the local cohomology modules.   
This is motivated by the work of the  second author with Verma in \cite{masuti-verma} for the filtration $\olin{\F}$.  In order to extend their result for the filtration $\F$, first we derive a formula for 
$[H^2_{(at_1,bt_2)}(\mathcal R^\prime(\F))]_{(r,s)}$
as  a direct limit of    $M^1_{r,s}(a^k,b^k)$ (Theorem~\ref{computation of local cohomology}).  
In \cite{masuti-verma} the authors prove that for the filtration $\olin{\F}$ and a good joint reduction $(a,b)$ of $\olin{\F}$
 $$[H^2_{(at_1,bt_2)}(\mathcal R^\prime(\olin{\F}))]_{(r,s)} \cong M^1_{r,s}(a^k,b^k;\olin{\F}):=\frac{\olin{I^{r+k}J^{s+k}}}{a^k\olin{I^rJ^{s+k}}+b^k \olin{I^{r+k} J^s}},$$ 
which in particular shows that the length of $[H^2_{(at_1,bt_2)}(\mathcal R^\prime(\olin{\F}))]_{(r,s)}$ is finite. However, for the filtration $\F$,  $[H^2_{(at_1,bt_2)}(\mathcal R^\prime(\F))]_{(r,s)}$ need not be equal to $M^1_{r,s}(a^k,b^k)$ for large $k$.  In fact, the length of $[H^2_{(at_1,bt_2)}(\mathcal R^\prime(\F))]_{(r,s)}$ can be infinite (see Example~\ref{example-not-finite}).  In this section we investigate the finiteness of the length of the local cohomology modules $[H^2_{(at_1,bt_2)}(\mathcal R^\prime(\F))]_{(r,s)}$ and their vanishing in terms of the Hilbert and bigraded Hilbert coefficients,  as well as in terms of the vanishing of the modules $M^1_{i,*}$ for all $i \geq r$ and $M^1_{*,j}$ for all $j \geq s$.

For the sake of simplicity we set $\mathcal{R}^\prime:=\mathcal{R}^\prime(\mathcal{F})$. Let  $a \in I$ and  $b\in J$.  Consider the Koszul co-complex
$$
K^{\bullet} ( (at_1)^k,(bt_2)^k; {\mathcal{R}^\prime}):
0 \lrar {\mathcal{R}^\prime}
\buildrel\alpha_k\over\lrar {\mathcal{R}^\prime}(k,0) \oplus
\mathcal{R}^\prime(0,k)
\buildrel\beta_k\over\lrar {\mathcal{R}^\prime}(k,k)
\lrar
0,
$$
where the maps are defined as,
$$
\alpha_k(1) = ((at_1)^k, (bt_2)^k)~~ ~~ \mbox{ and } \beta_k(u, v) =
-(bt_2)^ku + (at_1)^kv.
$$
\noindent 
Then for $i=0,1,2$,  
\begin{equation} \label{eq1}
H^i_{(at_1,bt_2)}({\mathcal{R}^\prime}) =
\displaystyle\lim_{\stackrel{\lrar}{k}}
H^i(K^{\bullet} ( (at_1)^k,(bt_2)^k; {\mathcal{R}^\prime})) \hspace{.2in} \mbox{\cite[Theorem~5.2.9]{bs}}.
\end{equation}

 In the 
following theorem we {recover} some results from \cite[Theorem 3.6]{masuti-verma}) for $[H^2_{(at_1,bt_2)}({{\mathcal{R}}^\prime})]_{(r,s)}$,  where $(a,b)$
is a Rees-joint reduction  of $I$ and $J$.

\begin{theorem}
 \label{computation of local cohomology} 
 Let $(R,\m)$ be a Cohen-Macaulay local ring of dimension two. Let  $I$ and $J$ be $\m$-primary 
 ideals in $R$ and $(a,b)$ a joint reduction of $I$ and $ J$.  Then for all 
$r,s \geq 0$,
\begin{enumerate}
 \item 
\label{computation of local cohomology-1} 
$
[H^2_{(at_1,bt_2)}({\mathcal{R}^\prime})]_{(r,s)} \cong
\displaystyle\lim_{\stackrel{\lrar}{k}}
  M^1_{r,s}(a^k,b^k).
$ 
 
\item 
\label{computation of local cohomology-2} 
If in addition $a$ and $b$ are Rees-superficial elements for $I$ and $J$,  respectively,  then for $k \gg 0$ the maps
$$
\mu_k:   M^1_{r,s}(a^k,b^k)
\buildrel{.(ab)}\over\lrar
 M^1_{r,s}(a^{k+1},b^{k+1})
$$ 
 are  injective. 
\end{enumerate}
\end{theorem}
\begin{proof}
(\ref{computation of local cohomology-1})  Follows from (\ref{eq1}).

(\ref{computation of local cohomology-2}) For $x \in I^{r+k}J^{s+k},$ we write $\olin{x}$ for the image of $x$ in $M^1_{r,s}(a^k,b^k)$.
Let $x \in I^{r+k}J^{s+k}$ be such that $\mu_k(\bar x)=0$. Then 
$xab=a^{k+1}p+b^{k+1}q$ for some $p \in I^rJ^{s+k+1}$ and $q\in I^{r+k+1}J^{s}. $ 
Hence $q \in (a) \cap I^{r+k+1}J^s=aI^{r+k}J^s$ for $k \gg 0$. Therefore 
$q=aq^\prime$ for some $q^\prime  \in I^{r+k}J^s$. Similarly for $k \gg 0$, $p=bp^\prime$ 
for some $p^\prime \in I^rJ^{s+k}$. Hence 
$$
x=a^kp^\prime+b^kq^\prime\in a^kI^rJ^{s+k}+b^kI^{r+k}J^s.
$$ 
Thus $\bar x=0$ and hence $\mu_k$ is injective for all $k \gg 0.$
\end{proof}

\noindent 
The map $\mu_k$ defined in Theorem~\ref{computation of local cohomology}  need not  be surjective  
for $k \gg0$. In fact, if  the map $\mu_k$ in Theorem~\ref{computation of local cohomology} is surjective 
for $k \gg 0$,  then $\mu_k$ is an isomorphism for $k \gg 0$ by Theorem~\ref{computation of local cohomology}  and  this implies that 
$[H^2_{(at_1,bt_2)}(\R^\prime)]_{(r,s)}$ has finite length which need not be true (see Example~\ref{example-not-finite}).   
The non-finiteness of the length of $[H^2_{(at_1,bt_2)}(\R^\prime)]_{(r,s)}$ is one of the obstructions in extending Rees' theorem for the ordinary powers of ideals. 
 In the following theorem we give  equivalent conditions for  $\ell_R([H^2_{(at_1,bt_2)}
 (\R^\prime)]_{(r,s)})$ to be finite in terms of the Hilbert and the bigraded Hilbert coefficients, and also in terms of the vanishing of the modules $M^1_{i,*}$ and $M^1_{*,j}.$
We remark that in \cite[Theorem 3.7]{masuti-verma} the authors derived a formula for $\ell([H^2_{(at_1,bt_2)}({{\mathcal{R}^\prime(\olin{\F})}})]_{(r,s)})$, which in particular shows that 
 $[H^2_{(at_1,bt_2)}({{\mathcal{R}^\prime(\olin{\F})}})]_{(r,s)}$ has finite length.

\begin{theorem}
\label{finite length criterion}
 Let $(R,\m)$ be a Cohen-Macaulay local ring of dimension two, $I$, $J$ be $\m$-primary ideals in $R$ and let $(a,b)$ a Rees-joint reduction of $I$ and $J$. 
Let $r_0,s_0 \geq0$ be fixed. Then  the following statements are equivalent:
\begin{enumerate}
\item \label{finite length criterion1}
$\ell_R([H^2_{(at_1,bt_2)}(\R^\prime)]_{(r_0,s_0)}) < \infty;$

\item \label{finite length criterion2}
$e_{(0,1)}=g_1(r_0)+r_0e_{(1,1)}$ and $e_{(1,0)}=h_1(s_0)+s_0e_{(1,1)}$; 

\item \label{finite length criterion3} $M^1_{i,*}=0$ for $i \geq r_0$ and $M^1_{*,j}=0$ for all $j \geq s_0.$
\end{enumerate}
If any of the above equivalent conditions hold true, then   
$$\ell_R([H^2_{(at_1,bt_2)}(\R^\prime)]_{(r_0,s_0)}) 
= -e_2(IJ)+g_2(r_0)+h_2(s_0)-\ell(R/I^{r_0}J^{s_0})+
r_0s_0e_{(1,1)}+\ell\left(\frac{\widetilde{\F}_{(a,b)}(I^{r_0}, J^{s_0})}{I^{r_0}J^{s_0}}\right).
$$

 \end{theorem}
\begin{proof}
(\ref{finite length criterion1}) $\implies$ (\ref{finite length criterion2}): 
By Proposition~\ref{computation of length}(\ref{k large}) for $k \gg 0$, $\ell(  M^1_{r_0,s_0}(a^k,b^k))$ is a polynomial in $k$
of degree at most $1$. 
By Theorem~\ref{computation of local cohomology} for all $k \gg 0$,
$$
\ell\left(M^1_{r_0,s_0}(a^k,b^k) \right) 
\leq \ell([H^2_{(at_1,bt_2)}({\R}^{\prime})]_{(r_0,s_0)}) < \infty.
$$
Hence  $\ell(M^1_{r_0,s_0}(a^k,b^k))$ is a constant for $k \gg0$. 
This implies that  
$[e_{(0,1)}-g_1(r_0)-r_0e_{(1,1)}]+[e_{(1,0)}-h_1(s_0)-s_0e_{(1,1)}]=0$ by Proposition~\ref{computation of length}(\ref{k large}). Since 
$e_{(0,1)}-g_1(r_0)-r_0e_{(1,1)}$ and $e_{(1,0)}-h_1(s_0)-s_0e_{(1,1)}$ are non-negative 
by Proposition~\ref{nonnegativity of coefficients}, 
$$
   e_{(0,1)} - g_1(r_0) - r_0e_{(1,1)}
= e_{(1,0)} - h_1(s_0) - s_0 e_{(1,1)}
= 0
$$
which gives the result. 

\noindent (\ref{finite length criterion2}) $\implies$ (\ref{finite length criterion3}): Follows from Proposition~\ref{nonnegativity of coefficients}.

\noindent (\ref{finite length criterion3}) $\implies$ (\ref{finite length criterion1}): 
Since $M^1_{i,*}=0$ for all $i \geq r_0$ and $M^1_{*,j}=0$ for all $j \geq s_0,$ by Proposition~\ref{nonnegativity of coefficients} 
\[
e_{(0,1)}=g_1(r_0)+r_0e_{(1,1)} \mbox{ and } e_{(1,0)}=h_1(s_0)+s_0 e_{(1,1)}.
\]
Substituting for $ e_{(0,1)}$ and $e_{(1,0)}$ in 
Proposition~\ref{computation of length}(\ref{k large}), we get $\ell( M^1_{r_0,s_0}(a^k,b^k) )$ is a constant for 
$k \gg 0$. 
Since $\mu_k$ is injective for $k \gg0,$ we conclude that 
$\mu_k$ is also surjective for $k \gg 0$ and hence is an isomorphism for $k \gg 0.$ Thus for $k \gg 0,$
$$
[H^2_{(at_1,bt_2)}(\R^{\prime})]_{(r_0,s_0)}
\cong
 M^1_{r_0,s_0}(a^k,b^k)  
$$
and hence $\ell([H^2_{(at_1,bt_2)}(\R^{\prime})]_{(r_0,s_0)}) < \infty$
\end{proof}

 \noindent As a corollary we give  equivalent conditions for 
$\ell_R([H^2_{(at_1,bt_2)}(\R^{\prime})]_{(0,0)})  <\infty$.  

\begin{corollary}
\label{corollary finite length criterion}
With the assumptions as in Theorem~\ref{finite length criterion}, 
the following statements are equivalent:
\begin{enumerate}
  \item   
 $\ell_R([H^2_{(at_1,bt_2)}(\R^{\prime})]_{(0,0)}) < \infty;$

\item 
$e_{(1,0)}=e_1(I)$ and $e_{(0,1)}=e_1(J);$

\item $M^1_{i,*}=0$ for $i \geq 0$ and $M^1_{*,j}=0$ for all $j \geq 0.$
 \end{enumerate}
 If any of the above equivalent conditions hold true, then    
$$\ell_R([H^2_{(at_1,bt_2)}(\R^{\prime})]_{(0,0)}) 
= -e_2(IJ)+e_2(I)+e_2(J).
$$
\end{corollary}
\begin{proof}
Put  $r_0=s_0=0$ in Theorem~\ref{finite length criterion}. 
\end{proof}

\begin{theorem}
\label{finiteness at any point implies finiteness at higher points}
With the assumptions as in Theorem~\ref{finite length criterion},  
if  
$\ell([H^2_{(at_1,bt_2)}(\R^{\prime})]_{(r_0,s_0)}) < \infty$ for some $r_0,s_0 \geq 0$, then 
$\ell([H^2_{(at_1,bt_2)}(\R^{\prime})]_{(r,s)}) < \infty$ for all  $r \geq r_0$ and $s \geq s_0$. 
\end{theorem}
\begin{proof}
Suppose $\ell([H^2_{(at_1,bt_2)}(\R^{\prime})]_{(r_0,s_0)})< \infty$ for some $r_0,s_0 \geq 0$. Then by  
Theorem~\ref{finite length criterion}(\ref{finite length criterion3}) $M^1_{i,*}=0$ for all $ i \geq r_0$ and $M^1_{*,j}=0$ for all $ j\geq s_0.$ 
As $r \geq r_0$ and $s \geq s_0,$ $M^1_{i,*}=0$ for all $ i \geq r$ and $M^1_{*,j}=0$ for all $ j\geq s.$ 
Therefore using Theorem~\ref{finite length criterion} once again we get 
$\ell([H^2_{(at_1,bt_2)}(\R^{\prime})]_{(r,s)}) < \infty$ for all  $r \geq r_0$ and $s \geq s_0$.
\end{proof}

\noindent In Section~\ref{Section:Examples} we will give an example to show that $e_{(1,0)} \neq e_1(I)$ and hence 
$\ell([H^2_{(at_1,bt_2)}(\R^{\prime})]_{(0,0)}) $ is not finite. Here, we give an example for which
$\ell_R([H^2_{(at_1,bt_2)}(\R^{\prime})]_{(r,s)})$ is finite, but need not be zero. Recall that an ideal $K \subseteq I$ is
called a {\it reduction} of $I$ if $K I^{n} = I^{n+1}$ for some $n$.

\begin{example}
\label{Example:finitelength}
 Let $(R,\m)$ be a Cohen-Macaulay local ring of dimension two and $I$ be an $\m$-primary 
ideal of $R$. Let $J=I$. Then $e_{(1,0)}=e_{(0,1)}=e_1(I)$. Hence, for any reduction $(a,b)$ of $I$ such that $a$ and $b$ are superficial elements, $\ell([H^2_{(at_1,bt_2)}(\R^{\prime})]_{(0,0)}) = e_2(I) <\infty$ 
 by Corollary~\ref{corollary finite length criterion}.
\end{example}

\noindent  In the rest of this section we give necessary and 
sufficient conditions for the vanishing of $[H^2_{(at_1,bt_2)}(\R^{\prime})]_{(r,s)}.$ 
The following theorem gives a cohomological interpretation of 
Theorem~\ref{arbitrary joint reduction number:KH}.

\begin{theorem}
\label{arbitrary joint reduction number}
Let $(R,\m)$ be a Cohen-Macaulay local ring of dimension two, $I,J$ be 
$\m$-primary ideals in $R$ and let $(a,b)$ a joint reduction of $I$ and $J$. Let  $r_0,s_0 \geq 0$.  Then the following statements are equivalent:
\begin{enumerate}
\item \label{arbitrary joint reduction number1} For all  $ k\gg 0$
$$I^{r_0+k}J^{s_0+k}=a^kI^{r_0}J^{s_0+k}+b^kI^{r_0+k}J^{s_0};$$

 \item \label{arbitrary joint reduction number2} $e_{(0,1)}=g_1(r_0) + r_0e_{(1,1)}$, $e_{(1,0)} = h_1(s_0) + s_0 e_{(1,1)}$ and 
 $e_2(IJ)=g_2(r_0)+h_2(s_0)-\ell(R/I^{r_0}J^{s_0})+
r_0s_0e_{(1,1)}+\ell\left(\frac{\widetilde{\F}_{(a,b)}(I^{r_0}, J^{s_0})}{I^{r_0}J^{s_0}}\right)$;

\item \label{arbitrary joint reduction number3}
  $[H^2_{(at_1,bt_2)}(\R^{\prime})]_{(r_0,s_0)}=0.$

\end{enumerate}
\end{theorem}
\begin{proof} 
The equivalence of (\ref{arbitrary joint reduction number1}) and (\ref{arbitrary joint reduction number2}) follows from Theorem~\ref{arbitrary joint reduction number:KH}, while the  equivalence of (\ref{arbitrary joint reduction number2}) and (\ref{arbitrary joint reduction number3}) follows from Theorem~\ref{finite length criterion}.
\end{proof}

\noindent
As a corollary we obtain a criteria for the vanishing of $[H^2_{(at_1,bt_2)}(\R^{\prime})]_{(0,0)}=0$. This  gives a 
cohomological interpretation of Theorem~\ref{theorem:computation of length}.

 \begin{corollary}
\label{joint reduction number zero for large powers}
Let $(R,\m)$ be a Cohen-Macaulay local ring of dimension two and $I,J$ be 
$\m$-primary ideals in $R$.  Then the following statements are equivalent:
\begin{enumerate}
\item 
$r(I^k| J^k)=0$ for $k \gg 0$;
 \item 
 $e_{(1,0)}=e_1(I)$, $e_{(0,1)}=e_1(J)$ and $e_2(IJ)=e_2(I)+e_2(J)$.
 
\item 
for every joint reduction $(a,b)$ of $I$ and $J$,
 $[H^2_{(at_1,bt_2)}(\R^{\prime})]_{(0,0)}=0$;
\end{enumerate}
\end{corollary}
\begin{proof}
 Put $r_0=s_0=0$ in Theorem~\ref{arbitrary joint reduction number}.
\end{proof}

In the next theorem we give a criteria for a vector to be in the joint reduction lattice of $I$ and $J$ in terms of the vanishing of 
$[H^2_{(at_1,bt_2)}(\R^{\prime})]_{(r_0,s_0)}$. This gives
a cohomological interpretation of Theorem~\ref{thm:ArbitraryJRN}. 

\begin{theorem}
\label{thm:ArbitraryJRN-LC}
 Let $(R,\m)$ be a Cohen-Macaulay local ring of dimension two and $I,J$ be 
$\m$-primary ideals in $R$.  Let  $r_0,s_0 \geq 0$. Assume that there exists a joint reduction  $(a,b)$ of $I$ and $J$ 
such that
\begin{eqnarray*}
 (a) \cap I^{r_0+k}J^{s_0}&=&aI^{r_0+k-1}J^{s_0} \mbox{ for }k \geq 1  \\
(b) \cap I^{r_0}J^{s_0+k}&=&bI^{r_0}J^{s_0+k-1} \mbox{ for }k \geq 1. 
\end{eqnarray*}
Then the following statements are equivalent:
\begin{enumerate}
 \item \label{thm:ArbitraryJRN1-LC1} $(r_0,s_0) \in \Lambda(I|J);$

 \item \label{thm:ArbitraryJRN2-LC2} $e_{(0,1)}=g_1(r_0) + r_0e_{(1,1)}$, $e_{(1,0)} = h_1(s_0) + s_0 e_{(1,1)}$ and 
 $e_2(IJ)=g_2(r_0)+h_2(s_0)-\ell(R/I^{r_0}J^{s_0})+
r_0s_0e_{(1,1)}+\ell\left(\frac{\widetilde{\F}_{(a,b)}(I^{r_0}, J^{s_0})}{I^{r_0}J^{s_0}}\right)$;

\item \label{thm:ArbitraryJRN3-LC3} $[H^2_{(at_1,bt_2)}(\R^{\prime})]_{(r_0,s_0)}=0.$
\end{enumerate}
\end{theorem}
\begin{proof}
The equivalence of (\ref{thm:ArbitraryJRN1-LC1}) and (\ref{thm:ArbitraryJRN2-LC2}) follows from Theorem~\ref{thm:ArbitraryJRN}, while the  equivalence of (\ref{thm:ArbitraryJRN2-LC2}) and (\ref{thm:ArbitraryJRN3-LC3}) follows from Theorem~\ref{arbitrary joint reduction number}.
\end{proof}

As a consequence we obtain the following corollary which gives a cohomological interpretation of Theorem~\ref{joint reduction number zero:KH}. 

\begin{corollary}
\label{joint reduction number zero}
  Let $(R,\m)$ be a Cohen-Macaulay local ring of dimension two and $I,J$ be 
$\m$-primary ideals in $R$. Assume that $\depth G(I),\depth G(J)\geq 1$. Then the following statements are equivalent.
\begin{enumerate}
\item \label{joint reduction number zero1}
$r(I|J) = 0$; 

 \item \label{joint reduction number zero2}
 $e_{(1,0)}=e_1(I)$, $e_{(0,1)}=e_1(J)$ and $e_2(IJ)=e_2(I)+e_2(J)$;

  \item \label{joint reduction number zero3}
for every  joint reduction $(a,b)$ of $I$ and $J$, $[H^2_{(at_1,bt_2)}(\R^{\prime})]_{(0,0)}=0$,
\end{enumerate}
\end{corollary}
\begin{proof}
By \cite[Lemma~1.2]{rees3} there exists a Rees-joint reduction 
 $(a,b)$ of $I$ and $J$ such that $a \in I\setminus \m I$ and $b \in J \setminus \m J$.  
 Let $a^*$ (resp. $b^*$) denotes the image of $a$ (resp. $b$) in $[G(I)]_1$ (resp. $[G(J)]_1$). 
 Since $\depth G(I)$ (resp. $\depth G(J)$) $\geq 1$, $a^*$ (resp. $b^*$) is a 
 nonzero-divisor in $G(I)$ (resp. $G(J)$)  by \cite[Lemma~2.1]{huckaba-marley}. 
 Hence 
 \begin{eqnarray*}
  (a) \cap I^n=aI^{n-1} \mbox{ and }(b) \cap J^n=bJ^{n-1} \mbox{ for all }n >0.
 \end{eqnarray*}
 Therefore the equivalence of (\ref{joint reduction number zero1}) and (\ref{joint reduction number zero2}) follows from Theorem~\ref{thm:ArbitraryJRN-LC} by taking $r_0=s_0=0.$ The equivalence of (\ref{joint reduction number zero2}) and (\ref{joint reduction number zero3}) follows from Corollary~\ref{joint reduction number zero for large powers}.
\end{proof}

In the following example we verify Theorem~\ref{joint reduction number zero:KH} and 
Corollary~\ref{joint reduction number zero} for complete ideals in a regular local ring.

\begin{example}
\label{Example:CIinRLR}
 Let $(R, \m)$ be a regular local ring of dimension two  and let $I,J$ be complete ideals 
(i.e. $\olin{I}=I$ and $\olin{J}=J$) in $R$. 
Then $IJ = aJ + bI$ for any joint reduction $(a,b)$ of $I$ and $J$ 
\cite[Theorem~2.1]{verma-joint}. 
Hence $M^1_{0,0}(a^k, b^k) = 0$ for all $k \geq 1$ which implies that 
$[H^2_{(at_1,bt_2)}(\R^{\prime})]_{(0,0)}=0$ by
Theorem~\ref{computation of local cohomology}(\ref{computation of local cohomology-1}). Moreover, by \cite[Theorem~3.2]{verma-joint} $e_{(1,0)} = e_1(I)$, $e_{(0,1)} = e_1(J)$ and $e_2(IJ)=e_2(I)+e_2(J)$.
This verifies Theorem~\ref{joint reduction number zero:KH} and 
Corollary~\ref{joint reduction number zero}.
\end{example}

\noindent 
\section{Examples}
\label{Section:Examples}
In this section we give an explicit example for which $e_{(0,1)} \neq e_1(J)$  and $[H^2_{(bt_1,at_2)}({\R^\prime})]_{(0,0)}$ is not finite  (Example~\ref{example-not-finite}). We also give an example where $e_{(1,0)} \neq e_1(I)$ and $e_{(0,1)} \neq e_1(J)$ (Example~\ref{Example:both}).

Recall that a reduction $K$
is called a {\it minimal reduction} of $I$ if whenever $K' \subseteq K$ and
$K'$ is a reduction of $I$, then $K'=K$  \cite{north-rees}. 
The {\it reduction number} of $I$ with respect to a
minimal reduction $K$ of $I$ is defined as 
$$ 
r_K (I) := \min\{ n \geq 0 \;|\; K I^n =I^{n+1} \}.  
$$ 
The {\it reduction number} of $I$ denoted by $r(I)$ is defined to be
the minimum of $r_K(I)$ where $K$ varies over all minimal
reductions of $I$. 

In order to obtain Example~\ref{example-not-finite} we need the following proposition.
\begin{proposition} 
\label{proposition for example}
 Let $(R, \m)$ be a Cohen-Macaulay local ring of dimension two and $I$ be an $\m$-primary ideal in 
 $R$ with $r(I) \geq 1$. Let $J=(a,b)$ be a minimal 
reduction of $I$ such that $b\in I$  (resp. $a\in J$) is a Rees-superficial element for $I$ (resp. $J$). 
Then  
\begin{enumerate}
 \item 
 \label{prop:example-a}
 $IJ^{s+1} \not = bJ^{s+1}+aIJ^{s}  \mbox{ for any }s \geq 0.$ 
 
 \item
  \label{prop:example-b}
  $e_{(0,1)} \neq e_1(J).$
\end{enumerate}
\end{proposition}
\begin{proof}
(\ref{prop:example-a})
Let  $I = (a,b, z_1, \ldots, z_{t-2})$ be a generating set of $I$ such that $z_i \notin (a,b) $ for every $i=1,\ldots,t-2.$.  
 Note that $t >2$ by \cite[Theorem~3.21]{marleythesis}. To prove the lemma it is enough to show that for all $s \geq 0$
\begin{eqnarray*}
b^{s+1} z_i \not \in bJ^{s+1}+aIJ^{s} \mbox{ for all }i=1, \ldots, t-2.
\end{eqnarray*}
Suppose $b^{s+1} z_i \in bJ^{s+1}+aIJ^{s}$ for some $i$. Inductively, for all $s \geq 0,$ we have
\begin{eqnarray} \nonumber
I J^{s+1}  &=& (a,b)^{s+2} + (a,b)^{s+1}  (z_1, \ldots,z_{t-2})  \mbox{ and }\\ \label{example-eqn1}
bJ^{s+1}+aIJ^{s} &=& (a,b)^{s+2} + a(a,b)^{s} (z_1, \ldots,z_{t-2}). 
\end{eqnarray}
 Hence from (\ref{example-eqn1}), 
\begin{eqnarray*}
b^{s+1} z_i 
= \sum_{k=0}^{s+2} x_k a^k b^{s+2-k}  
+ a \sum_{k=0}^{s} \left( \sum_{j=1}^{t-2}  y_{kj}a^kb^{s-k}z_j\right)
\end{eqnarray*}
where $x_k,y_{kj} \in R$. This implies that 
\begin{eqnarray*}
b^{s+1}(z_i  - bx_0) \in (a).
\end{eqnarray*}
As $a,b$ is a regular sequence in $R$, 
 \begin{eqnarray*}
 z_i - bx_0\in (a).
 \end{eqnarray*}
 Therefore $z_i \in (a,b)$ which contradicts that $z_i \notin (a,b).$

(\ref{prop:example-b})
  Suppose $e_{(0,1)}=e_1(J)$. As $(b,a)$ is a Rees-joint reduction of $I$ and $J$, by Corollary~\ref{inequality between e's}(\ref{inequality between e's-2}) $M^1_{i,s}(b,a)$ for all $i \geq 0$ and $s \gg 0$.  In particular, 
  \beqn
  0= M^1_{0,s}(b,a)
= \frac{I J^{s+1}}
                       {b J^{s+1}+ a IJ^{s}} . 
                       \eeqn
This contradicts (\ref{prop:example-a}).
 \end{proof}
 
\noindent Now we give an explicit example for which $e_{(0,1)}\neq e_1(J)$. In fact, the following example shows that the difference $e_{(0,1)}-e_1(J)$ can be as large as possible. 
 
 \begin{example}
 \label{example-not-finite}
  Let $R=k[\!|x,y|\!]$, $\m= (x,y)$, $I=\m^t$, $J=(x^t,y^t)$, $t \geq 2$. Put $a= x^t$ and $b = y^t$. 
Then $b \in I$ (resp. $a \in J$) is a Rees-superficial element for $I$ and $J$.  
  Therefore by Proposition~\ref{proposition for example}
  $e_{(0,1)} \neq e_1(J)$. We explicitly calculate  $e_{(0,1)} -e_1(J)$.
 For all $r, s \geq 1$, 
  \beqn
  \ell \left(  \f{R} {I^r J^s}\right) 
 & =&  \ell \left(  \f{R} {\m^{t(r+s)}}\right)\\
&  =& {t(r+s) + 1 \choose 2}\\
  &= &t^2{r+1 \choose 2} + t^2 rs + t^2{s+1 \choose 2}
  - {t \choose 2}r - {t \choose 2} s.
  \eeqn
  As $J$ is a parameter ideal $e_1(J)=0$. Hence 
  $e_{(0,1)}- e_1(J) = {t \choose 2}$. 
Therefore by Corollary~\ref{corollary finite length criterion} the length of  
    $[H^2_{(bt_1,at_2)}({\R^\prime})]_{(0,0)}$ is not finite.
  \qed
  \end{example}

Notice that in Example~\ref{example-not-finite} $e_{(1,0)}=e_1(I).$ Next we give an example for which $e_{(1,0)} \neq e_1(I)$ as well as $e_{(0,1)} \neq e_1(J).$

\begin{example}
\label{Example:both}
Let $R = k[\!|x,y|\!]_{\m}$ where $\m = (x,y)$ and let $a \geq 1$. Put $I = (x^2, y^2)$ and $J = (x^{3}, y^{3})$. We claim that 
\been
\item
\label{Example:both-one}
$e_{(1,0)} \not = e_1(I)$, $e_{(0,1)} \not = e_1(J)$, $e_2(IJ) \not = e_2(I) + e_2(J)$.

\item
\label{Example:both-two}
$r(I |J) \not =0$ even though $G(I)$ and $G(J)$ are Cohen-Macaulay.

\item
\label{Example:both-three}
$[H^2_{( (x^2)t_1,(y^3)t_2)}({\mathcal{R}^\prime})]_{(0,0)} \neq 0.$
\eeen
(\ref{Example:both-one}) One can verify that 
$
(IJ)^{3}
= (x^2 J + y^2 J)^{3}
\subseteq x^2 I^2J^{3} +y^{3} I^{3}J^2
\subseteq  (IJ)^{3}.
$
Hence  $I^3J^3= x^2 I^2J^{3} +y^{3} I^{3}J^2$
and  for all $r,s \geq 3$,
\beq
\label{joint-reduction}
I^r J^s 
= x^2 I^{r-1}J^s + y^3 I^r J^{s-1}
=x^{2(r-2)} I^{2}J^s + y^{3(s-2)} I^r J^{2}.
\eeq
 For the rest of the proof we will assume that $r,s \geq 3$. Recall   the complex $C_{\bullet}^{\prime} ((x^{2(r-2)}, y^{3(s-2)}),2,2)$ from (\ref{Eqn:Complex2})
\beq
\label{example:complex1}
       0
\lrar \frac{R}{I^2 J^2}
\stackrel{\psi_1 }{\lrar}  
   \f{R}{I^r J^2}  \oplus   \frac{R}{I^2 J^s} 
\stackrel{\psi_0 }{\lrar} 
      \frac{ (x^{2(r-2)} , y^{3(s-2) })} {x^{2(r-2)} I^{2}J^s + y^{3(s-2)} I^r J^{2}}
\lrar 0
\eeq
where ${\displaystyle \psi_1:=\small \left(\begin{array}{c} x^{2(r-2)}\\ y^{3(s-2)}\end{array}\right)}$ and 
$\psi_0:=\small\left(\begin{array}{cc}y^{3(s-2)}&-x^{2(r-2)} \end{array}\right)$.
\normalsize
We claim that the complex (\ref{example:complex1}) is exact. 
Clearly $\psi_0$ is surjective. From the proof of Theorem~\ref{computation of homology}(\ref{computation of homology c}),  $\ker(\psi_0) = \im (\psi_1).$ We show that $\psi_1$ is injective.
From the complex (\ref{example:complex1}),   $\ker(\psi_1) := (I^2 J^s : y^{3(s-2)} ) \cap (I^r J^2 : x^{2(r-2)}) / I^2 J^2 $. Therefore for $s \geq 3$
\beq
\label{colon-1}
(x^9y) + I^2J^2   \nonumber
&\subseteq &  (I^2 J^s : y^{3(s-2)} )\\ \nonumber
&=&   (x^9 I^2 J^{s-3} + y^{3(s-2)} I^2 J^2 : y^{3(s-2)}   )\\ \nonumber
&=& (x^9 ( x^2 I + (y^4)) J^{(s-3)} + y^{3(s-2)} I^2 J^2 : y^{3(s-2)} ) \\ \nonumber
&=& \begin{cases}
(x^{11})  + (x^9y) + I^2 J^2 & \mbox{ if } s=3\\
(x^{11})  + (x^{12}, x^9y) + I^2 J^2 & \mbox{ if } s>3\\
\end{cases}\\
 &\subseteq&  (x^9y) +  I^2 J^2.  
\eeq
Hence  the equality holds in (\ref{colon-1}). 
Let $r \geq 3$. Then
\beq
\label{colon-2} 
(xy^9) + I^2J^2
&\subseteq&      (I^r J^2 : x^{2(r-2)})\\ \nonumber
&=& (x^{2(r-2)} I^2 J^2 +y^6 I^{2(r-3)} J^2:x^{2(r-2)}) \\ \nonumber
&=& \begin{cases}
I^2 J^2 + y^6(x^4, xy^3, y^6)  & \mbox{ if } r=3\\
I^2 J^2 + y^6(x^4,  x^2y^2, xy^3, y^5) & \mbox{ if } r=4\\
I^2 J^2 + y^6(x^4,  x^2y^2, xy^3, y^4) & \mbox{ if } r>4\\
\end{cases}\\  
&\subseteq& (xy^9) + I^2J^2.
\eeq
Hence  the equality holds in (\ref{colon-2}). 
From (\ref{colon-1}) and (\ref{colon-2}) we get 
\beqn
 (I^2 J^s : y^{3(s-2)} ) \cap (I^r J^2 : x^{2(r-2)})
 = ( (xy^9) + I^2 J^2) \cap ( (x^9y) + I^2 J^2)
 = (x^9y^9) + I^2 J^2 = I^2 J^2. 
\eeqn
This implies that $\psi_1$ is injective. From (\ref{joint-reduction}) and (\ref{example:complex1}), for all $r,s \geq 3$, 
\beqn
\ell \left(  \frac{R} {I^rJ^s}\right) 
= \ell \left(  \frac{R}  {(x^{2(r-2)} , y^{3(s-2)} )} \right)
+     \ell \left( \frac{R}{I^2 J^s} \right) 
+   \ell \left(\frac{R}{I^r J^2} \right)
   -    \ell \left(\frac{R}{I^2 J^2} \right)
\eeqn
By induction  we can show that for all $r,s \geq 2$,
\beqn
I^{2} J^s &=& (x^{4 + 3s}, y^{4 + 3s}) + x^2 y^2 (x,y)^{3s} \mbox{ and }\\
I^r J^2 &= & (x^{2r+6}, y^{2r + 6}) + x^2y^2(x,y)^{2r+2}.
\eeqn
Therefore 
\beqn
        \ell \left ( \f{R}{I^2J^s} \right)
&=& \ell \left( \f{R}{\m^{4+3s}}\right)+2 
=       \binom{4+3s+1}{2} + 2\mbox{ and } \\
         \ell \left( \f{R}{I^rJ^2} \right) 
&=& \ell \left( \f{R}{\m^{2r+6}} \right)+2 
         = \binom{2r+6+1}{2}+2.
\eeqn
Hence for all $r,s \geq 3$,
\beq
\label{Eqn:LenghComputation}
\ell \left(  \frac{R} {I^rJ^s}\right) 
&=& 6(r-2) (s-2)
 +  {4+3s +1 \choose 2} + 2
+ {2r+6 +1 \choose 2} + 2
 - {11 \choose 2} - 2\\ \nonumber 
  &=& 4 {r+1 \choose 2} + 6rs +  9 {s+1 \choose 2} - ( r + 3s) +2.
\eeq
Hence $e_{(1,0)} = 1$, $e_{(0,1)} = 3$. 
As $I$ and $J$ are parameter ideals $e_1(I) =e_1(J)  = e_2(I) = e_2(J)= 0$. 
Therefore $e_{(1,0)}  - e_1(I)= 1$,  $e_{(0,1)} - e_1(J)= 3$ and $e_2(IJ) -e_2(I) -e_2(J) = 2$.

(\ref{Example:both-two})
By Theorem~\ref{joint reduction number zero:KH}  $r(I|J) \not = 0$,  even though $G(I)$ and $G(J)$ are Cohen-Macaulay.
We verify this directly here. 
 Suppose  joint reduction number of $I$ and $J$ is zero. Then by \cite[Theorem~3.2(d)]{verma-joint}, 
$e_{(1,1)}(I|J) = \ell(R/ IJ) - \ell (R/I) - \ell (R/J)$.  
But here
$$
   \ell \left( \f{R}{ IJ} \right) - \ell \left(\f{R}{I} \right) - \ell \left( \f{R}{J} \right) 
   = 17-4-9 
= 4 \neq 6 
=  e_{(1,1)}(I|J).$$

 (\ref{Example:both-three}) Applying   (\ref{Example:both-one}) and (\ref{Example:both-two}) to Theorem~\ref{joint reduction number zero}  we conclude that $[H^2_{( (x^2)t_1,(y^3)t_2)}({\mathcal{R}^\prime})]_{(0,0)} \neq 0.$ We verify this directly.  By Theorem~\ref{computation of local cohomology}(\ref{computation of local cohomology-1})
\beqn
[H^2_{( (x^2)t_1,(y^3)t_2)}({\mathcal{R}^\prime})]_{(0,0)} \cong
\displaystyle\lim_{\stackrel{\lrar}{k}}
  M^1_{0,0}((x^2)^k,(y^3)^k).
  \eeqn
 For all $ k \geq 3$, 
  \beq \nonumber
&&\ell \left( M^1_{0,0}((x^2)^k,(y^3)^k) \right)\\\nonumber
&=& \ell \left( \frac{I^{k} J^{k}}
                       {(x^2)^k J^{k}+ (y^3)^k I^{k}} \right)\\ \nonumber
&=& \ell \left( \frac{R} {(x^{2k}, y^{3k})}\right)+\ell\left(\frac{(x^{2k}, y^{3k})}{x^{2k}J^k+y^{3k}I^k}\right) - \ell\left( \frac{R}{I^kJ^k} \right)\\ \nonumber
&=& \ell \left( \frac{R} {(x^{2k}, y^{3k})}\right) 
+  \ell \left( \frac{R} {I^{k} }\right)
+ \ell \left( \frac{R} {J^{k}}\right)
-  \ell \left( \frac{R} {I^{k} J^{k}}\right)  \hspace*{1cm} \mbox{[by \cite[Lemma 3.1]{verma-joint}]}\\      \nonumber
&=&   6k^2 + 4\binom{k+1}{2} + 9 \binom{k+1}{2} 
        - \left[  4 {k+1 \choose 2} + 6k^2 +  9 {k+1 \choose 2} - 4k  + 2\right] \hspace*{1cm} \mbox{[From (\ref{Eqn:LenghComputation})]}\\ \nonumber
&=& 4k-2.
                       \eeq
As $k \geq 3$,   $4k-2 \not =0$. Hence  $[H^2_{( (x^2)t_1,(y^3)t_2)}({\mathcal{R}^\prime})]_{(0,0)} \neq 0$.
           \end{example}

\end{document}